\documentclass{amsart}
\usepackage{amssymb}

\usepackage{amsmath}
\usepackage{amscd}
\usepackage{thmdefs}

\input{tcilatex}
\theoremstyle{definition}
\theoremstyle{remark}
\numberwithin{equation}{section}

\input tcilatex

\begin{document}
\title[Graph $W^{*}$-Probability Spaces]{Graph $W^{*}$-Probability Spaces}
\author{Ilwoo Cho}
\address{Univ. of Iowa, Dep. of Math, Iowa City, IA, U. S. A}
\email{ilcho@math.uiowa.edu}
\date{}
\dedicatory{}
\thanks{}
\maketitle

\begin{abstract}
In this paper, we constructe a $W^{*}$-probability space $\left(
W^{*}(G),E\right) $ with amalgamation over a von Neumann algebra $D_{G},$
where $W^{*}(G)$ is a graph $W^{*}$-algebra induced by the countable
directed graph $G$. In this structure, we compute the $D_{G}$-valued moments
and cumulants of arbitrary random variables, by using the lattice path
models and we characterize the $D_{G}$-freeness of generators of $W^{*}(G),$
by the so-called diagram-distinctness on $G.$ As examples, we will compute
the $D_{G}$-moments and $D_{G}$-cumulants of the generating operators $%
T_{1}^{N}$ of $W^{*}(G_{1}^{N})$ and $T_{N}$ of $W^{*}(C_{N}),$ where $%
G_{1}^{N}$ is the graph with one vertex and $N$-distinct loop-edges and $%
C_{N}$ is the circulant graph with $N$-vertices and $N$-edges. In
particular, the generating operator $T_{1}^{N}$ is the free sum of $N$%
-semicircular elements. This operator $T_{1}^{N}$ is semicircular with its
even moments $\left( 2N\right) ^{k}\cdot c_{k},$ for all $k\in 2\Bbb{N},$
and its second cumulant $2N,$ where $c_{k}$ is the $k$-th Catalan number.
\end{abstract}

\strut

In this paper, we construct the graph $W^{*}$-probability spaces. The graph $%
W^{*}$-probability theory is one of the good example of Speicher's
combinatorial free probability theory with amalgamation (See [16]). In this
paper, we will observe how to compute the moment and cumulant of an
arbitrary random variables in the graph $W^{*}$-probability space and the
freeness on it with respect to the given conditional expectation. In [10],
Kribs and Power defined the free semigroupoid algebras and obtained some
properties of them. Our work is highly motivated by [10]. Roughly speaking,
graph $W^{*}$-algebras are $W^{*}$-topology closed version of free
semigroupoid algebras. Throughout this paper, let $G$ be a countable
directed graph and let $\mathbb{F}^{+}(G)$ be the free semigroupoid of $G,$
in the sense of Kribs and Power. i.e., it is a collection of all vertices of
the graph $G$ as units and all admissible finite paths, under the
admissibility. As a set, the free semigroupoid $\mathbb{F}^{+}(G)$ can be
decomposed by

\strut

\begin{center}
$\mathbb{F}^{+}(G)=V(G)\cup FP(G),$
\end{center}

\strut

where $V(G)$ is the vertex set of the graph $G$ and $FP(G)$ is the set of
all admissible finite paths. Trivially the edge set $E(G)$ of the graph $G$
is properly contained in $FP(G),$ since all edges of the graph can be
regarded as finite paths with their length $1.$ We define a graph $W^{*}$%
-algebra of $G$ by

\strut

\begin{center}
$W^{*}(G)\overset{def}{=}\overline{%
\mathbb{C}[\{L_{w},L_{w}^{*}:w\in
\mathbb{F}^{+}(G)\}]}^{w},$
\end{center}

\strut

where $L_{w}$ and $L_{w}^{\ast }$ are creation operators and annihilation
operators on the generalized Fock space $H_{G}=l^{2}\left( \mathbb{F}%
^{+}(G)\right) $ induced by the given graph $G,$ respectively. Notice that
the creation operators induced by vertices are projections and the creation
operators induced by finite paths are partial isometries. We can define the $%
W^{\ast }$-subalgebra $D_{G}$ of $W^{\ast }(G),$ which is called the
diagonal subalgebra by

\strut

\begin{center}
$D_{G}\overset{def}{=}\overline{\mathbb{C}[\{L_{v}:v\in V(G)\}]}^{w}.$
\end{center}

\strut

Then each element $a$ in the graph $W^{*}$-algebra $W^{*}(G)$ is expressed by

\strut

\begin{center}
$a=\underset{w\in \mathbb{F}^{+}(G:a),\,u_{w}\in \{1,*\}}{\sum }%
p_{w}L_{w}^{u_{w}},$ \ for $p_{w}\in \mathbb{C},$
\end{center}

\strut

where $\mathbb{F}^{+}(G:a)$ is a support of the element $a$, as a subset of
the free semigroupoid $\mathbb{F}^{+}(G).$ The above expression of the
random variable $a$ is said to be the Fourier expansion of $a.$ Since $%
\mathbb{F}^{+}(G)$ is decomposed by the disjoint subsets $V(G)$ and $FP(G),$
the support $\mathbb{F}^{+}(G:a)$ of $a$ is also decomposed by the following
disjoint subsets,

\strut

\begin{center}
$V(G:a)=\mathbb{F}^{+}(G:a)\cap V(G)$
\end{center}

and

\begin{center}
$FP(G:a)=\mathbb{F}^{+}(G:a)\cap FP(G).$
\end{center}

\strut

Thus the operator $a$ can be re-expressed by

\strut

\begin{center}
$a=\underset{v\in V(G:a)}{\sum }p_{v}L_{v}+\underset{w\in FP(G:a),\,u_{w}\in
\{1,*\}}{\sum }p_{w}L_{w}^{u_{w}}.$
\end{center}

\strut

Notice that if $V(G:a)\neq \emptyset ,$ then $\underset{v\in V(G:a)}{\sum }%
p_{v}L_{v}$ is contained in the diagonal subalgebra $D_{G}.$ Thus we have
the canonical conditional expectation $E:W^{*}(G)\rightarrow D_{G},$ defined
by

\strut

\begin{center}
$E\left( a\right) =\underset{v\in V(G:a)}{\sum }p_{v}L_{v},$
\end{center}

\strut

for all $a=\underset{w\in \mathbb{F}^{+}(G:a),\,u_{w}\in \{1,*\}}{\sum }%
p_{w}L_{w}^{u_{w}}$ \ in $W^{*}(G).$ Then the algebraic pair $\left(
W^{*}(G),E\right) $ is a $W^{*}$-probability space with amalgamation over $%
D_{G}$ (See [16]). It is easy to check that the conditional expectation $E$
is faithful in the sense that if $E(a^{*}a)=0_{D_{G}},$ for $a\in W^{*}(G),$
then $a=0_{D_{G}}.$

\strut

For the fixed operator $a\in W^{*}(G),$ the support $\mathbb{F}^{+}(G:a)$ of
the operator $a$ is again decomposed by

\strut

\begin{center}
$\mathbb{F}^{+}(G:a)=V(G:a)\cup FP_{*}(G:a)\cup FP_{*}^{c}(G:a),$
\end{center}

\strut

with the decomposition of $FP(G:a),$

$\strut $

\begin{center}
$FP(G:a)=FP_{*}(G:a)\cup FP_{*}^{c}(G:a),$
\end{center}

where

\strut

\begin{center}
$FP_{*}(G:a)=\{w\in FP(G:a):$both $L_{w}$ and $L_{w}^{*}$ are summands of $%
a\}$
\end{center}

and

\begin{center}
$FP_{*}(G:a)=FP(G:a)\,\,\setminus \,\,FP_{*}(G:a).$
\end{center}

\strut

The above new expression plays a key role to find the $D_{G}$-valued moments
of the random variable $a.$ In fact, the summands $p_{v}L_{v}$'s and $%
p_{w}L_{w}+p_{w^{t}}L_{w}^{*},$ for $v\in V(G:a)$ and $w\in FP_{*}(G:a)$ act
for the computation of $D_{G}$-valued moments of $a.$ By using the above
partition of the support of a random variable, we can compute the $D_{G}$%
-valued moments and $D_{G}$-valued cumulants of it via the lattice path
model $LP_{n}$ and the lattice path model $LP_{n}^{*}$ satisfying the $*$%
-axis-property. At a first glance, the computations of $D_{G}$-valued
moments and cumulants look so abstract (See Chapter 3) and hence it looks
useless. However, these computations, in particular the computation of $%
D_{G} $-valued cumulants, provides us how to figure out the $D_{G}$-freeness
of random variables by making us compute the mixed cumulants. As
applications, in the final chapter, we can compute the moment and cumulant
of the operator that is the sum of $N$-free semicircular elements with their
covariance $2.$ If $a$ is the operator, then the $n$-th moment of $a$ is

\strut

\begin{center}
$\left\{ 
\begin{array}{ll}
\left( 2N\right) ^{\frac{n}{2}}\cdot c_{\frac{n}{2}} & \text{if }n\text{ is
even} \\ 
0 & \text{if }n\text{ is odd,}
\end{array}
\right. $
\end{center}

\strut

and the $n$-th cumulant of $a$ is

\strut

\begin{center}
$\left\{ 
\begin{array}{lll}
2N &  & \text{if }n=2 \\ 
0 &  & \text{otherwise,}
\end{array}
\right. $
\end{center}

\strut \strut \strut

where $c_{k}=\frac{1}{k+1}\left( 
\begin{array}{l}
2k \\ 
\,\,k
\end{array}
\right) $ is the $k$-th Catalan number.

\strut

\strut Based on the $D_{G}$-cumulant computation, we can characterize the $%
D_{G}$-freeness of generators of $W^{*}(G),$ by the so-called
diagram-distinctness on the graph $G.$ i.e., the random variables $L_{w_{1}}$
and $L_{w_{2}}$ are free over $D_{G}$ if and only if $w_{1}$ and $w_{2}$ are
diagram-distinct the sense that $w_{1}$ and $w_{2}$ have different diagrams
on the graph $G.$ Also, we could find the necessary condition for the $D_{G}$%
-freeness of two arbitrary random variables $a$ and $b.$ i.e., if the
supports $\mathbb{F}^{+}(G:a)$ and $\mathbb{F}^{+}(G:b)$ are
diagram-distinct, in the sense that $w_{1}$ and $w_{2}$ are diagram distinct
for all pairs $(w_{1},w_{2})$ $\in $ $\mathbb{F}^{+}(G:a)$ $\times $ $%
\mathbb{F}^{+}(G:b),$ then the random variables $a$ and $b$ are free over $%
D_{G}.$

\strut \strut

\strut \strut

\strut \strut

\section{Graph $W^{*}$-Probability Spaces}

\strut

\strut

Let $G$ be a countable directed graph and let $\Bbb{F}^{+}(G)$ be the free
semigroupoid of $G.$ i.e., the set $\mathbb{F}^{+}(G)$ is the collection of
all vertices as units and all admissible finite paths of $G.$ Let $w$ be a
finite path with its source $s(w)=x$ and its range $r(w)=y,$ where $x,y\in
V(G).$ Then sometimes we will denote $w$ by $w=xwy$ to express the source
and the range of $w.$ We can define the graph Hilbert space $H_{G}$ by the
Hilbert space $l^{2}\left( \mathbb{F}^{+}(G)\right) $ generated by the
elements in the free semigroupoid $\mathbb{F}^{+}(G).$ i.e., this Hilbert
space has its Hilbert basis $\mathcal{B}=\{\xi _{w}:w\in \mathbb{F}%
^{+}(G)\}. $ Suppose that $w=e_{1}...e_{k}\in FP(G)$ is a finite path with $%
e_{1},...,e_{k}\in E(G).$ Then we can regard $\xi _{w}$ as $\xi
_{e_{1}}\otimes ...\otimes \xi _{e_{k}}.$ So, in [10], Kribs and Power
called this graph Hilbert space the generalized Fock space. Throughout this
paper, we will call $H_{G}$ the graph Hilbert space to emphasize that this
Hilbert space is induced by the graph.

\strut

Define the creation operator $L_{w},$ for $w\in \mathbb{F}^{+}(G),$ by the
multiplication operator by $\xi _{w}$ on $H_{G}.$ Then the creation operator 
$L$ on $H_{G}$ satisfies that

\strut

(i) \ $L_{w}=L_{xwy}=L_{x}L_{w}L_{y},$ for $w=xwy$ with $x,y\in V(G).$

\strut

(ii) $L_{w_{1}}L_{w_{2}}=\left\{ 
\begin{array}{lll}
L_{w_{1}w_{2}} &  & \text{if }w_{1}w_{2}\in \mathbb{F}^{+}(G) \\ 
&  &  \\ 
0 &  & \text{if }w_{1}w_{2}\notin \mathbb{F}^{+}(G),
\end{array}
\right. $

\strut

\ \ \ \ for all $w_{1},w_{2}\in \mathbb{F}^{+}(G).$

\strut

Now, define the annihilation operator $L_{w}^{*},$ for $w\in \mathbb{F}%
^{+}(G)$ by

\strut

\begin{center}
$L_{w}^{\ast }\xi _{w^{\prime }}\overset{def}{=}\left\{ 
\begin{array}{lll}
\xi _{h} &  & \text{if }w^{\prime }=wh\in \mathbb{F}^{+}(G)\xi \\ 
&  &  \\ 
0 &  & \text{otherwise.}
\end{array}
\right. $
\end{center}

\strut

The above definition is gotten by the following observation ;

\strut

\begin{center}
$
\begin{array}{ll}
<L_{w}\xi _{h},\xi _{wh}>\, & =\,<\xi _{wh},\xi _{wh}>\, \\ 
& =\,1=\,<\xi _{h},\xi _{h}> \\ 
& =\,<\xi _{h},L_{w}^{*}\xi _{wh}>,
\end{array}
\,$
\end{center}

\strut

where $<,>$ is the inner product on the graph Hilbert space $H_{G}.$ Of
course, in the above formula we need the admissibility of $w$ and $h$ in $%
\mathbb{F}^{+}(G).$ However, even though $w$ and $h$ are not admissible
(i.e., $wh\notin \mathbb{F}^{+}(G)$), by the definition of $L_{w}^{\ast },$
we have that

\strut

\begin{center}
$
\begin{array}{ll}
<L_{w}\xi _{h},\xi _{h}> & =\,<0,\xi _{h}> \\ 
& =0=\,<\xi _{h},0> \\ 
& =\,<\xi _{h},L_{w}^{*}\xi _{h}>.
\end{array}
\,\,$
\end{center}

\strut

Notice that the creation operator $L$ and the annihilation operator $L^{*}$
satisfy that

\strut

(1.1) \ \ \ $L_{w}^{*}L_{w}=L_{y}$ \ \ and \ \ $L_{w}L_{w}^{*}=L_{x},$ \ for
all \ $w=xwy\in \mathbb{F}^{+}(G),$

\strut

\textbf{under the weak topology}, where $x,y\in V(G).$ Remark that if we
consider the von Neumann algebra $W^{*}(\{L_{w}\})$ generated by $L_{w}$ and 
$L_{w}^{*}$ in $B(H_{G}),$ then the projections $L_{y}$ and $L_{x}$ are
Murray-von Neumann equivalent, because there exists a partial isometry $%
L_{w} $ satisfying the relation (1.1). Indeed, if $w=xwy$ in $\mathbb{F}%
^{+}(G), $ with $x,y\in V(G),$ then under the weak topology we have that

\strut

(1,2) \ \ \ $L_{w}L_{w}^{*}L_{w}=L_{w}$ \ \ and \ \ $%
L_{w}^{*}L_{w}L_{w}^{*}=L_{w}^{*}.$

\strut

So, the creation operator $L_{w}$ is a partial isometry in $W^{*}(\{L_{w}\})$
in $B(H_{G}).$ Assume now that $v\in V(G).$ Then we can regard $v$ as $%
v=vvv. $ So,

\strut

(1.3) $\ \ \ \ \ \ \ \ \ L_{v}^{*}L_{v}=L_{v}=L_{v}L_{v}^{*}=L_{v}^{*}.$

\strut

This relation shows that $L_{v}$ is a projection in $B(H_{G})$ for all $v\in
V(G).$

\strut

Define the \textbf{graph }$W^{*}$\textbf{-algebra} $W^{*}(G)$ by

\strut

\begin{center}
$W^{*}(G)\overset{def}{=}\overline{%
\mathbb{C}[\{L_{w},L_{w}^{*}:w\in
\mathbb{F}^{+}(G)\}]}^{w}.$
\end{center}

\strut

Then all generators are either partial isometries or projections, by (1.2)
and (1.3). So, this graph $W^{\ast }$-algebra contains a rich structure, as
a von Neumann algebra. (This construction can be the generalization of that
of group von Neumann algebra.) Naturally, we can define a von Neumann
subalgebra $D_{G}\subset W^{\ast }(G)$ generated by all projections $L_{v},$ 
$v\in V(G).$ i.e.

\strut

\begin{center}
$D_{G}\overset{def}{=}W^{*}\left( \{L_{v}:v\in V(G)\}\right) .$
\end{center}

\strut

We call this subalgebra the \textbf{diagonal subalgebra} of $W^{*}(G).$
Notice that $D_{G}=\Delta _{\left| G\right| }\subset M_{\left| G\right| }(%
\mathbb{C}),$ where $\Delta _{\left| G\right| }$ is the subalgebra of $%
M_{\left| G\right| }(\mathbb{C})$ generated by all diagonal matrices. Also,
notice that $1_{D_{G}}=\underset{v\in V(G)}{\sum }L_{v}=1_{W^{*}(G)}.$

\strut

If $a\in W^{*}(G)$ is an operator, then it has the following decomposition
which is called the Fourier expansion of $a$ ;

\strut

(1.4) $\ \ \ \ \ \ \ \ \ \ \ a=\underset{w\in \mathbb{F}^{+}(G:a),\,u_{w}\in
\{1,*\}}{\sum }p_{w}L_{w}^{u_{w}},$

\strut

where $p_{w}\in C$ and $\mathbb{F}^{+}(G:a)$ is the support of $a$ defined by

\strut

\begin{center}
$\mathbb{F}^{+}(G:a)=\{w\in \mathbb{F}^{+}(G):p_{w}\neq 0\}.$
\end{center}

\strut

Remark that the free semigroupoid $\mathbb{F}^{+}(G)$ has its partition $%
\{V(G),FP(G)\},$ as a set. i.e.,

\strut

\begin{center}
$\mathbb{F}^{+}(G)=V(G)\cup FP(G)$ \ \ and \ \ $V(G)\cap FP(G)=\emptyset .$
\end{center}

\strut

So, the support of $a$ is also partitioned by

\strut

\begin{center}
$\mathbb{F}^{+}(G:a)=V(G:a)\cup FP(G:a),$
\end{center}

\strut where

\begin{center}
$V(G:a)\overset{def}{=}V(G)\cap \mathbb{F}^{+}(G:a)$
\end{center}

and

\begin{center}
$FP(G:a)\overset{def}{=}FP(G)\cap \mathbb{F}^{+}(G:a).$
\end{center}

\strut

So, the above Fourier expansion (1.4) of the random variable $a$ can be
re-expressed by

\strut

(1.5) $\ \ \ \ \ \ a=\underset{v\in V(G:a)}{\sum }p_{v}L_{v}+\underset{w\in
FP(G:a),\,u_{w}\in \{1,*\}}{\sum }p_{w}L_{w}^{u_{w}}.$

\strut

We can easily see that if $V(G:a)\neq \emptyset ,$ then $\underset{v\in
V(G:a)}{\sum }p_{v}L_{v}$ is contained in the diagonal subalgebra $D_{G}.$
Also, if $V(G:a)=\emptyset ,$ then $\underset{v\in V(G:a)}{\sum }%
p_{v}L_{v}=0_{D_{G}}.$ So, we can define the following canonical conditional
expectation $E:W^{*}(G)\rightarrow D_{G}$ by

\strut

(1.6) \ \ \ $E(a)=E\left( \underset{w\in \mathbb{F}^{+}(G:a),\,u_{w}\in
\{1,*\}}{\sum }p_{w}L_{w}^{u_{w}}\right) \overset{def}{=}\underset{v\in
V(G:a)}{\sum }p_{v}L_{v},$

\strut

for all $a\in W^{*}(G).$ Indeed, $E$ is a well-determined conditional
expectation ; it is a bimodule map satisfying that

\strut

\begin{center}
$E(d)=d,$ for all $d\in D_{G}.$
\end{center}

\strut

And

\strut

\begin{center}
$
\begin{array}{ll}
E\left( dad^{\prime }\right) & =E\left( d(a_{d}+a_{0})d^{\prime }\right)
=E\left( da_{d}d^{\prime }+da_{0}d^{\prime }\right) \\ 
& =E\left( da_{d}d^{\prime }\right) =da_{d}d^{\prime }=d\left( E(a)\right)
d^{\prime },
\end{array}
$
\end{center}

\strut

for all $d,d^{\prime }\in D_{G}$ and $a=a_{d}+a_{0}\in W^{*}(G),$ where

\strut

\begin{center}
$a_{d}=\underset{v\in V(G:a)}{\sum }p_{v}L_{v}$ \ \ and \ \ $a_{0}=\underset{%
w\in FP(G:a),\,u_{w}\in \{1,*\}}{\sum }p_{w}L_{w}^{u_{w}}.$
\end{center}

\strut

Also,

\strut

\begin{center}
$E\left( a^{\ast }\right) =E\left( (a_{d}+a_{0})^{\ast }\right) =E\left(
a_{d}^{\ast }+a_{0}^{\ast }\right) =a_{d}^{\ast }=E(a)^{\ast },$
\end{center}

\strut

for all $a\in W^{*}(G).$ Here, $a_{d}^{*}=\left( \underset{v\in V(G:a)}{\sum 
}p_{v}L_{v}\right) ^{*}=\underset{v\in V(G:a)}{\sum }\overline{p_{v}}\,L_{v}$
in $D_{G}.$

\strut \strut \strut \strut

\begin{definition}
Let $G$ be a countable directed graph and let $W^{*}(G)$ be the graph $W^{*}$%
-algebra induced by $G.$ Let $E:W^{*}(G)\rightarrow D_{G}$ be the
conditional expectation defined above. Then we say that the algebraic pair $%
\left( W^{*}(G),E\right) $ is the graph $W^{*}$-probability space over the
diagonal subalgebra $D_{G}$. By the very definition, it is one of the $W^{*}$%
-probability space with amalgamation over $D_{G}.$ All elements in $\left(
W^{*}(G),E\right) $ are called $D_{G}$-valued random variables.
\end{definition}

\strut

We have a graph $W^{*}$-probability space $\left( W^{*}(G),E\right) $ over
its diagonal subalgebra $D_{G}.$ We will define the following free
probability data of $D_{G}$-valued random variables.

\strut

\begin{definition}
Let $W^{*}(G)$ be the graph $W^{*}$-algebra induced by $G$ and let $a\in
W^{*}(G).$ Define the $n$-th ($D_{G}$-valued) moment of $a$ by

\strut

$\ \ \ \ \ E\left( d_{1}ad_{2}a...d_{n}a\right) ,$ for all $n\in \mathbb{N}$,

\strut

where $d_{1},...,d_{n}\in D_{G}$. Also, define the $n$-th ($D_{G}$-valued)
cumulant of $a$ by

\strut

$\ \ \ \ \ k_{n}(d_{1}a,d_{2}a,...,d_{n}a)=C^{(n)}\left( d_{1}a\otimes
d_{2}a\otimes ...\otimes d_{n}a\right) ,$

\strut

for all $n\in \mathbb{N},$ and for $d_{1},...,d_{n}\in D_{G},$ where $%
\widehat{C}=(C^{(n)})_{n=1}^{\infty }\in I^{c}\left( W^{*}(G),D_{G}\right) $
is the cumulant multiplicative bimodule map induced by the conditional
expectation $E,$ in the sense of Speicher. We define the $n$-th trivial
moment of $a$ and the $n$-th trivial cumulant of $a$ by

\strut

$\ \ \ \ \ E(a^{n})$ $\ \ $and $\ \ k_{n}\left( \underset{n-times}{%
\underbrace{a,a,...,a}}\right) =C^{(n)}\left( a\otimes a\otimes ...\otimes
a\right) ,$

\strut

respectively, for all $n\in \mathbb{N}.$
\end{definition}

\strut

To compute the $D_{G}$-valued moments and cumulants of the $D_{G}$-valued
random variable $a,$ we need to introduce the following new definition ;

\strut

\begin{definition}
Let $\left( W^{*}(G),E\right) $ be a graph $W^{*}$-probability space over $%
D_{G}$ and let $a\in \left( W^{*}(G),E\right) $ be a random variable. Define
the subset $FP_{*}(G:a)$ in $FP(G:a)$ \ by

\strut

$\ \ \ FP_{*}\left( G:a\right) \overset{def}{=}\{w\in \mathbb{F}^{+}(G:a):$%
both $L_{w}$ and $L_{w}^{*}$ are summands of $a\}.$

\strut

And let $FP_{*}^{c}(G:a)\overset{def}{=}FP(G:a)\,\setminus \,FP_{*}(G:a).$
\end{definition}

\strut \strut \strut

We already observed that if $a\in \left( W^{*}(G),E\right) $ is a $D_{G}$%
-valued random variable, then $a$ has its Fourier expansion $a_{d}+a_{0},$
where

\strut

\begin{center}
$a_{d}=\underset{v\in V(G:a)}{\sum }p_{v}L_{v}$
\end{center}

and

\begin{center}
$a_{0}=\underset{w\in FP(G:a),\,u_{w}\in \{1,*\}}{\sum }p_{w}L_{w}^{u_{w}}.$
\end{center}

\strut

By the previous definition, the set $FP(G:a)$ is partitioned by

\strut

\begin{center}
$FP(G:a)=FP_{*}(G:a)\cup FP_{*}^{c}(G:a),$
\end{center}

\strut

for the fixed random variable $a$ in $\left( W^{*}(G),E\right) .$ So, the
summand $a_{0},$ in the Fourier expansion of $a=a_{d}+a_{0},$ has the
following decomposition ;

\strut

\begin{center}
$a_{0}=a_{(*)}+a_{(non-*)},$
\end{center}

\strut where\strut

\begin{center}
$a_{(*)}=\underset{l\in FP_{*}(G:a)}{\sum }\left(
p_{l}L_{l}+p_{l^{t}}L_{l}^{*}\right) $
\end{center}

and

\begin{center}
$a_{(non-*)}=\underset{w\in FP_{*}^{c}(G:a),\,u_{w}\in \{1,*\}}{\sum }%
p_{w}L_{w}^{u_{w}},$
\end{center}

\strut

where $p_{l^{t}}$ is the coefficient of $L_{l}^{\ast }$ depending on $l\in
FP_{\ast }(G:a).$ (There is no special meaning for the complex number $%
p_{l^{t}}.$ But we have to keep in mind that $p_{l}\neq p_{l^{t}},$ in
general. i.e. $a_{(\ast )}=\underset{l_{1}\in FP_{\ast }(G:a)}{\sum }%
p_{l_{1}}L_{l_{1}}+\underset{l_{2}\in FP_{\ast }(G:a)}{\sum }%
p_{l_{2}}L_{l_{2}}^{\ast }$ ! But for the convenience of using notation, we
will use the notation $p_{l^{t}},$ for the coefficient of $L_{l}^{\ast }.$)
For instance, let $V(G:a)=\{v_{1},v_{2}\}$ and $FP(G:a)=\{w_{1},w_{2}\}$ and
let the random variable $a$ in $\left( W^{\ast }(G),E\right) $ be

\strut

\begin{center}
$a=L_{v_{1}}+L_{v_{2}}+L_{w_{1}}^{*}+L_{w_{1}}+L_{w_{2}}^{*}.$
\end{center}

\strut

Then we have that $a_{d}=L_{v_{1}}+L_{v_{2}}$, $%
a_{(*)}=L_{w_{1}}^{*}+L_{w_{1}}$ and $a_{(non-*)}=L_{w_{2}}^{*}.$ By
definition, $a_{0}=a_{(*)}+a_{(non-*)}.$

\strut \strut \strut

\strut \strut

\strut \strut

\section{$D_{G}$-Moments and $D_{G}$-Cumulants of Random Variables}

\strut

\strut

\strut

Throughout this chapter, let $G$ be a countable directed graph and let $%
\left( W^{*}(G),E\right) $ be the graph $W^{*}$-probability space over its
diagonal subalgebra $D_{G}.$ In this chapter, we will compute the $D_{G}$%
-valued moments and the $D_{G}$-valued cumulants of arbitrary random variable

$\strut $

\begin{center}
$a=\underset{w\in \mathbb{F}^{+}(G:a),\,u_{w}\in \{1,*\}}{\sum }%
p_{w}L_{w}^{u_{w}}$
\end{center}

\strut

in the graph $W^{*}$-probability space $\left( W^{*}(G),E\right) $.

\strut

\strut

\strut

\subsection{Lattice Path Model}

\strut

\strut

\strut

Throughout this section, let $G$ be a countable directed graph and let $%
\left( W^{*}(G),E\right) $ be the graph $W^{*}$-probability space over its
diagonal subalgebra $D_{G}.$ Let $w_{1},...,w_{n}\in \Bbb{F}^{+}(G)$ and let 
$L_{w_{1}}^{u_{w_{1}}}...L_{w_{n}}^{u_{w_{n}}}\in \left( W^{*}(G),E\right) $
be a $D_{G}$-valued random variable. In this section, we will define a
lattice path model for the random variable $%
L_{w_{1}}^{u_{w_{1}}}...L_{w_{n}}^{u_{w_{n}}}.$ Recall that if $%
w=e_{1}....e_{k}\in FP(G)$ with $e_{1},...,e_{k}\in E(G),$ then we can
define the length $\left| w\right| $ of $w$ by $k.$ i.e.e, the length $%
\left| w\right| $ of $w$ is the cardinality $k$ of the admissible edges $%
e_{1},...,e_{k}.$

\strut

\begin{definition}
Let $G$ be a countable directed graph and $\Bbb{F}^{+}(G),$ the free
semigroupoid. If $w\in \Bbb{F}^{+}(G),$ then $L_{w}$ is the corresponding $%
D_{G}$-valued random variable in $\left( W^{*}(G),E\right) .$ We define the
lattice path $l_{w}$ of $L_{w}$ and the lattice path $l_{w}^{-1}$ of $%
L_{w}^{*}$ by the lattice paths satisfying that ;

\strut

(i) \ \ the lattice path $l_{w}$ starts from $*=(0,0)$ on the $\Bbb{R}^{2}$%
-plane.

\strut

(ii) \ if $w\in V(G),$ then $l_{w}$ has its end point $(0,1).$

\strut

(iii) if $w\in E(G),$ then $l_{w}$ has its end point $(1,1).$

\strut

(iv) if $w\in E(G),$ then $l_{w}^{-1}$ has its end point $(-1,-1).$

\strut

(v) \ if $w\in FP(G)$ with $\left| w\right| =k,$ then $l_{w}$ has its end
point $(k,k).$

\strut

(vi) if $w\in FP(G)$ with $\left| w\right| =k,$ then $l_{w}^{-1}$ has its
end point $(-k,-k).$

\strut

Assume that finite paths $w_{1},...,w_{s}$ in $FP(G)$ satisfy that $%
w_{1}...w_{s}\in FP(G).$ Define the lattice path $l_{w_{1}...w_{s}}$ by the
connected lattice path of the lattice paths $l_{w_{1}},$ ..., $l_{w_{s}}.$
i.e.e, $l_{w_{2}}$ starts from $(k_{w_{1}},k_{w_{1}})\in \Bbb{R}^{+}$ and
ends at $(k_{w_{1}}+k_{w_{2}},k_{w_{1}}+k_{w_{2}}),$ where $\left|
w_{1}\right| =k_{w_{1}}$ and $\left| w_{2}\right| =k_{w_{2}}.$ Similarly, we
can define the lattice path $l_{w_{1}...w_{s}}^{-1}$ as the connected path
of $l_{w_{s}}^{-1},$ $l_{w_{s-1}}^{-1},$ ..., $l_{w_{1}}^{-1}.$
\end{definition}

\strut

\begin{definition}
Let $G$ be a countable directed graph and assume that $%
L_{w_{1}},...,L_{w_{n}}$ are generators of $\left( W^{*}(G),E\right) .$ Then
we have the lattice paths $l_{w_{1}},$ ..., $l_{w_{n}}$ of $L_{w_{1}},$ ..., 
$L_{w_{n}},$ respectively in $\Bbb{R}^{2}.$ Suppose that $%
L_{w_{1}}^{u_{w_{1}}}...L_{w_{n}}^{u_{w_{n}}}\neq 0_{D_{G}}$ in $\left(
W^{*}(G),E\right) ,$ where $u_{w_{1}},...,u_{w_{n}}\in \{1,*\}.$ Define the
lattice path $l_{w_{1},...,w_{n}}^{u_{w_{1}},...,u_{w_{n}}}$ of nonzero $%
L_{w_{1}}^{u_{w_{1}}}...L_{w_{n}}^{u_{w_{n}}}$ by the connected lattice path
of $l_{w_{1}}^{t_{w_{1}}},$ ..., $l_{w_{n}}^{t_{w_{n}}},$ where $t_{w_{j}}=1$
if $u_{w_{j}}=1$ and $t_{w_{j}}=-1$ if $u_{w_{j}}=*.$ Assume that $%
L_{w_{1}}^{u_{w_{1}}}...L_{w_{n}}^{u_{w_{n}}}$ $=$ $0_{D_{G}}.$ Then the
empty set $\emptyset $ in $\Bbb{R}^{2}$ is the lattice path of it. We call
it the empty lattice path. By $LP_{n},$ we will denote the set of all
lattice paths of the $D_{G}$-valued random variables having their forms of $%
L_{w_{1}}^{u_{w_{1}}}...L_{w_{n}}^{u_{w_{n}}},$ including empty lattice path.
\end{definition}

\strut

Also, we will define the following important property on the set of all
lattice paths ;

\strut

\begin{definition}
Let $l_{w_{1},...,w_{n}}^{u_{w_{1}},...,u_{w_{n}}}\neq \emptyset $ be a
lattice path of $L_{w_{1}}^{u_{w_{1}}}...L_{w_{n}}^{u_{w_{n}}}\neq 0_{D_{G}}$
in $LP_{n}.$ If the lattice path $%
l_{w_{1},...,w_{n}}^{u_{w_{1}},...,u_{w_{n}}}$ starts from $*$ and ends on
the $*$-axis in $\Bbb{R}^{+},$ then we say that the lattice path $%
l_{w_{1},...,w_{n}}^{u_{w_{1}},...,u_{w_{n}}}$ has the $*$-axis-property. By 
$LP_{n}^{*},$ we will denote the set of all lattice paths having their forms
of $l_{w_{1},...,w_{n}}^{u_{w_{1}},...,u_{w_{n}}}$ which have the $*$%
-axis-property. By little abuse of notation, sometimes, we will say that the 
$D_{G}$-valued random variable $%
L_{w_{1}}^{u_{w_{1}}}...L_{w_{n}}^{u_{w_{n}}} $satisfies the $*$%
-axis-property if the lattice path $%
l_{w_{1},...,w_{n}}^{u_{w_{1}},...,u_{w_{n}}}$ of it has the $*$%
-axis-property.
\end{definition}

\strut

The following theorem shows that finding $E\left(
L_{w_{1}}^{u_{w_{1}}}...L_{w_{n}}^{u_{w_{n}}}\right) $ is checking the $*$%
-axis-property of $L_{w_{1}}^{u_{w_{1}}}...L_{w_{n}}^{u_{w_{n}}}.$

\strut

\begin{theorem}
Let $L_{w_{1}}^{u_{w_{1}}}...L_{w_{n}}^{u_{w_{n}}}\in \left(
W^{*}(G),E\right) $ be a $D_{G}$-valued random variable, where $%
u_{w_{1}},...,u_{w_{n}}\in \{1,*\}.$ Then $E\left(
L_{w_{1}}^{u_{w_{1}}}...L_{w_{n}}^{u_{w_{n}}}\right) $ $\neq $ $0_{D_{G}}$
if and only if $L_{w_{1}}^{u_{w_{1}}}...L_{w_{n}}^{u_{w_{n}}}$ has the $*$%
-axis-property (i.e., the corresponding lattice path $%
l_{w_{1},...,w_{n}}^{u_{w_{1}},...,u_{w_{n}}}$ of $%
L_{w_{1}}^{u_{w_{1}}}...L_{w_{n}}^{u_{w_{n}}}$ is contained in $LP_{n}^{*}.$
Notice that $\emptyset \notin LP_{n}^{*}.$)
\end{theorem}

\strut

\begin{proof}
($\Leftarrow $) Let $l=l_{w_{1},...,w_{n}}^{u_{w_{1}},...,u_{w_{n}}}\in
LP_{n}^{*}.$ Suppose that $w_{1}=vw_{1}v_{1}^{\prime }$ and $%
w_{n}=v_{n}w_{n}v_{n}^{\prime },$ for $v_{1},$ $v_{1}^{\prime },$ $v_{n},$ $%
v_{n}^{\prime }$ $\in $ $V(G).$ If $l$ is in $LP_{n}^{*},$ then

\strut

(2.1.1)$\ \ \ \ \ \ \ \left\{ 
\begin{array}{lll}
v_{1}=v_{n}^{\prime } &  & \text{if }u_{w_{1}}=1\text{ and }u_{w_{n}}=1 \\ 
&  &  \\ 
v_{1}=v_{n} &  & \text{if }u_{w_{1}}=1\text{ and }u_{w_{n}}=* \\ 
&  &  \\ 
v_{1}^{\prime }=v_{n}^{\prime } &  & \text{if }u_{w_{1}}=*\text{ and }%
u_{w_{n}}=1 \\ 
&  &  \\ 
v_{1}^{\prime }=v_{n} &  & \text{if }u_{w_{1}}=*\text{ and }u_{w_{n}}=*.
\end{array}
\right. $

\strut

By the definition of lattice paths having the $*$-axis-property and by
(2.1.1), if $l_{w_{1},...,w_{n}}^{u_{w_{1}},...,u_{w_{n}}}\in LP_{n}^{*},$
then there exists $v\in V(G)$ such that

\strut

$\ \ \ \ \ \ \ \ \ \ \ \ \ \ \ \ \ \ \
L_{w_{1}}^{u_{w_{1}}}...L_{w_{n}}^{u_{w_{n}}}=L_{v},$

where

(2.1.2) \ \ \ \ $\ \ \left\{ 
\begin{array}{ll}
v=v_{1}=v_{n}^{\prime } & \text{if }u_{w_{1}}=1\text{ and }u_{w_{n}}=1 \\ 
&  \\ 
v=v_{1}=v_{n} & \text{if }u_{w_{1}}=1\text{ and }u_{w_{n}}=* \\ 
&  \\ 
v=v_{1}^{\prime }=v_{n}^{\prime } & \text{if }u_{w_{1}}=*\text{ and }%
u_{w_{n}}=1 \\ 
&  \\ 
v=v_{1}^{\prime }=v_{n} & \text{if }u_{w_{1}}=*\text{ and }u_{w_{n}}=*.
\end{array}
\right. $

\strut

This shows that $E\left(
L_{w_{1}}^{u_{w_{1}}}...L_{w_{n}}^{u_{w_{n}}}\right) =L_{v}\neq 0_{D_{G}}.$

\strut

($\Rightarrow $) Assume that $E\left(
L_{w_{1}}^{u_{w_{1}}}...L_{w_{n}}^{u_{w_{n}}}\right) \neq 0_{D_{G}}.$ This
means that there exists $L_{v},$ with $v\in V(G),$ such that

\strut

(2.1.3) \ \ \ \ \ \ \ \ $\ \ \ E\left(
L_{w_{1}}^{u_{w_{1}}}...L_{w_{n}}^{u_{w_{n}}}\right) =L_{v}.$

\strut

Equivalently, we have that $%
L_{w_{1}}^{u_{w_{1}}}...L_{w_{n}}^{u_{w_{n}}}=L_{v}$ in $W^{*}(G).$ Let $%
l=l_{w_{1},...,w_{n}}^{u_{w_{1}},...,u_{w_{n}}}\in LP_{n}$ be the lattice
path of the $D_{G}$-valued random variable $%
L_{w_{1}}^{u_{w_{1}}}...L_{w_{n}}^{u_{w_{n}}}.$ By (2.1.3), trivially, $%
l\neq \emptyset ,$ since $l$ should be the connected lattice path. Assume
that the nonempty lattice path $l$ is contained in $LP_{n}\,\,\setminus
\,LP_{n}^{*}.$ Then, under the same conditions of (2.1.1), we have that

\strut

(2.1.4) \ \ \ \ \ \ $\left\{ 
\begin{array}{lll}
v_{1}\neq v_{n}^{\prime } &  & \text{if }u_{w_{1}}=1\text{ and }u_{w_{n}}=1
\\ 
&  &  \\ 
v_{1}\neq v_{n} &  & \text{if }u_{w_{1}}=1\text{ and }u_{w_{n}}=* \\ 
&  &  \\ 
v_{1}^{\prime }\neq v_{n}^{\prime } &  & \text{if }u_{w_{1}}=*\text{ and }%
u_{w_{n}}=1 \\ 
&  &  \\ 
v_{1}^{\prime }\neq v_{n} &  & \text{if }u_{w_{1}}=*\text{ and }u_{w_{n}}=*.
\end{array}
\right. $

\strut

Therefore, by (2.1.2), there is no vertex $v$ satisfying $%
L_{w_{1}}^{u_{w_{1}}}...L_{w_{n}}^{u_{w_{n}}}=L_{v}.$ This contradict our
assumption.
\end{proof}

\strut

By the previous theorem, we can conclude that $E\left(
L_{w_{1}}^{u_{w_{1}}}...L_{w_{n}}^{u_{w_{n}}}\right) =L_{v},$ for some $v\in
V(G)$ if and only if the lattice path $%
l_{w_{1},...,w_{n}}^{u_{w_{1}},...,u_{w_{n}}}$ has the $*$-axis-property
(i.e., $l_{w_{1},...,w_{n}}^{u_{w_{1}},...,u_{w_{n}}}\in LP_{n}^{*}$).\strut
\strut \strut

\strut \strut

\subsection{$D_{G}$-Valued Moments and Cumulants of Random Variables\strut}

\bigskip

\bigskip

Let $w_{1},...,w_{n}\in \mathbb{F}^{+}(G)$, $u_{1},...,u_{n}\in \{1,*\}$ and
let $L_{w_{1}}^{u_{1}}...L_{w_{n}}^{u_{n}}\in \left( W^{*}(G),E\right) $ be
a $D_{G}$-valued random variable. Recall that, in the previous section, we
observed that the $D_{G}$-valued random variable $%
L_{w_{1}}^{u_{1}}...L_{w_{n}}^{u_{n}}=L_{v}\in \left( W^{*}(G),E\right) $
with $v\in V(G)$ if and only if the lattice path $%
l_{w_{1},...,w_{n}}^{u_{1},...,u_{n}}$ of \ $%
L_{w_{1}}^{u_{1}}...L_{w_{n}}^{u_{n}}$ has the $*$-axis-property
(equivalently, $l_{w_{1},...,w_{n}}^{u_{1},...,u_{n}}\in LP_{n}^{*}$).
Throughout this section, fix a $D_{G}$-valued random variable $a\in \left(
W^{*}(G),E\right) .$ Then the $D_{G}$-valued random variable $a$ has the
following Fourier expansion,

\bigskip

\begin{center}
$a=\underset{v\in V(G:a)}{\sum }p_{v}L_{v}+\underset{l\in FP_{*}(G:a)}{\sum }%
\left( p_{l}L_{l}+p_{l^{t}}L_{l}\right) +\underset{w\in
FP_{*}^{c}(G:a),~u_{w}\in \{1,*\}}{\sum }p_{w}L_{w}^{u_{w}}.$
\end{center}

\bigskip

Let's observe the new $D_{G}$-valued random variable $d_{1}ad_{2}a...d_{n}a%
\in \left( W^{*}(G),E\right) ,$ where $d_{1},...,d_{n}\in D_{G}$ and $a\in
W^{*}(G)$ is given. Put

\strut

\begin{center}
$d_{j}=\underset{v_{j}\in V(G:d_{j})}{\sum }q_{v_{j}}L_{v_{j}}\in D_{G},$
for \ $j=1,...,n.$
\end{center}

\strut

Notice that $V(G:d_{j})=\mathbb{F}^{+}(G:d_{j}),$ since $d_{j}\in
D_{G}\hookrightarrow W^{\ast }(G).$ Then

\strut

$\ d_{1}ad_{2}a...d_{n}a$

\strut

$\ \ \ =\left( \underset{v_{1}\in V(G:d_{1})}{\sum }q_{v_{1}}L_{v_{1}}%
\right) \left( \underset{w_{1}\in \mathbb{F}^{+}(G:a),\,u_{w_{1}}\in \{1,*\}%
}{\sum }p_{w_{1}}L_{w_{1}}^{u_{w_{1}}}\right) $

$\ \ \ \ \ \ \ \ \ \ \ \ \ \ \ \ \ \cdot \cdot \cdot \left( \underset{%
v_{1}\in V(G:d_{n})}{\sum }q_{v_{n}}L_{v_{n}}\right) \left( \underset{%
w_{n}\in \mathbb{F}^{+}(G:a),\,u_{w_{n}}\in \{1,*\}}{\sum }%
p_{w_{n}}L_{w_{n}}^{u_{w_{n}}}\right) $

\strut

$\ \ \ =\underset{(v_{1},...,v_{n})\in \Pi _{j=1}^{n}V(G:d_{j})}{\sum }%
\left( q_{v_{1}}...q_{v_{n}}\right) $

$\ \ \ \ \ \ \ \ \ \ \ \ \ \ \ \ \ \ \ \ \ \ \ \ \ \ \ \ \ \ \ \
(L_{v_{1}}\left( \underset{w_{1}\in \mathbb{F}^{+}(G:a),\,u_{w_{1}}\in
\{1,*\}}{\sum }p_{w_{1}}L_{w_{1}}^{u_{w_{1}}}\right) $

$\ \ \ \ \ \ \ \ \ \ \ \ \ \ \ \ \ \ \ \ \ \ \ \ \ \ \ \ \ \ \ \ \ \ \ \ \ \
\ \cdot \cdot \cdot L_{v_{n}}\left( \underset{w_{n}\in \mathbb{F}%
^{+}(G:a),\,u_{w_{n}}\in \{1,*\}}{\sum }p_{w_{n}}L_{w_{n}}^{u_{w_{n}}}%
\right) )$

\strut

(1.2.1)

\strut

$\ \ \ =\underset{(v_{1},...,v_{n})\in \Pi _{j=1}^{n}V(G:d_{j})}{\sum }%
\left( q_{v_{1}}...q_{v_{n}}\right) $

\strut

$\ \ \ \ \ \underset{(w_{1},...,w_{n})\in \mathbb{F}^{+}(G:a)^{n},\,%
\,u_{w_{j}}\in \{1,*\}}{\sum }\left( p_{w_{1}}...p_{w_{n}}\right)
L_{v_{1}}L_{w_{1}}^{u_{w_{1}}}...L_{v_{n}}L_{w_{n}}^{u_{w_{n}}}.$

\strut

Now, consider the random variable $%
L_{v_{1}}L_{w_{1}}^{u_{w_{1}}}...L_{v_{n}}L_{w_{n}}^{u_{w_{n}}}$ in the
formula (1.2.1). Suppose that $w_{j}=x_{j}w_{j}y_{j},$ with $x_{j},y_{j}\in
V(G),$ for all \ $j=1,...,n.$ Then

\strut

\strut (1.2.2)

\begin{center}
$L_{v_{1}}L_{w_{1}}^{u_{w_{1}}}...L_{v_{n}}L_{w_{n}}^{u_{w_{n}}}=\delta
_{(v_{1},x_{1},y_{1}:u_{w_{1}})}\cdot \cdot \cdot \delta
_{(v_{n},x_{n},y_{n}:u_{w_{n}})}\left(
L_{w_{1}}^{u_{w_{n}}}...L_{w_{n}}^{u_{w_{n}}}\right) ,$
\end{center}

\strut

where

\strut

\begin{center}
$\delta _{(v_{j},x_{j},y_{j}:u_{w_{j}})}=\left\{ 
\begin{array}{lll}
\delta _{v_{j},x_{j}} &  & \text{if }u_{w_{j}}=1 \\ 
&  &  \\ 
\delta _{v_{j},y_{j}} &  & \text{if }u_{w_{j}}=*,
\end{array}
\right. $
\end{center}

\strut

for all \ $j=1,...,n,$ where $\delta $ in the right-hand side is the
Kronecker delta. So, the left-hand side can be understood as a (conditional)
Kronecker delta depending on $\{1,*\}$.

\strut \strut

By (1.2.1) and (1.2.2), the $n$-th moment of $a$ is

\strut

$\ E\left( d_{1}a...d_{n}a\right) $

\strut

\ \ $\ \ \ =E(\underset{(v_{1},...,v_{n})\in \Pi _{j=1}^{n}V(G:d_{j})}{\sum }%
\left( \Pi _{j=1}^{n}q_{v_{j}}\right) $

\strut

$\ \ \ \ \ \ \ \ \ \ \ \ \ \ \ \underset{(w_{1},...,w_{n})\in \mathbb{F}%
^{+}(G:a)^{n},\,w_{j}=x_{j}w_{j}y_{j},\,u_{w_{j}}\in \{1,\ast \}}{\sum }%
\left( \Pi _{j=1}^{n}p_{w_{j}}\right) $

\strut

$\ \ \ \ \ \ \ \ \ \ \ \ \ \ \ \ \ \ \ \ \ \ \ \ \ \ \ \ \ \ \ \ \ \ \left(
\Pi _{j=1}^{n}\delta _{(v_{j},x_{j},y_{j}:u_{w_{j}})}\right) \left(
L_{w_{1}}^{u_{w_{1}}}...L_{w_{n}}^{u_{w_{n}}}\right) )$

\strut

$\ \ \ \ \ \ =\underset{(v_{1},...,v_{n})\in \Pi _{j=1}^{n}V(G:d_{j})}{\sum }%
\left( \Pi _{j=1}^{n}q_{v_{j}}\right) $

\strut

$\ \ \ \ \ \ \ \ \ \ \underset{(w_{1},...,w_{n})\in \mathbb{F}%
^{+}(G:a)^{n},\,w_{j}=x_{j}w_{j}y_{j},\,u_{w_{j}}\in \{1,\ast \}}{\sum }%
\left( \Pi _{j=1}^{n}p_{w_{j}}\right) $

\strut

$\ \ \ \ \ \ \ \ \ \ \ \ \ \ \ \ \ \ \ \ \ \ \ \ \ \ \ \ \ \ \ \ \left( \Pi
_{j=1}^{n}\delta _{(v_{j},x_{j},y_{j}:u_{w_{j}})}\right) \,\,E\left(
L_{w_{1}}^{u_{w_{1}}}...L_{w_{n}}^{u_{w_{n}}}\right) .$

\strut

Thus to compute the $n$-th moment of $a$, we have to observe $E\left(
L_{w_{1}}^{u_{w_{1}}}...L_{w_{n}}^{u_{w_{n}}}\right) .$ In the previous
section, we observed that $E\left(
L_{w_{1}}^{u_{w_{1}}}...L_{w_{n}}^{u_{w_{n}}}\right) $ is nonvanishing if
and only if $L_{w_{1}}^{u_{w_{1}}}...L_{w_{n}}^{u_{w_{n}}}$ has the $*$%
-axis-property.

\strut

\begin{proposition}
Let $a\in \left( W^{*}(G),E\right) $ be given as above. Then the $n$-th
moment of $a$ is

\strut

$\ \ E\left( d_{1}a...d_{n}a\right) =\underset{(v_{1},...,v_{n})\in \Pi
_{j=1}^{n}V(G:d_{j})}{\sum }\left( \Pi _{j=1}^{n}q_{v_{j}}\right) $

\strut

$\ \ \ \underset{(w_{1},...,w_{n})\in \mathbb{F}^{+}(G:a)^{n},\,u_{w_{j}}\in
\{1,*\},\,l_{w_{1},...,w_{n}}^{u_{w_{1}},...,u_{w_{n}}}\in LP_{n}^{*}}{\sum }%
\left( \Pi _{j=1}^{n}p_{w_{j}}\right) $

\strut

$\ \ \ \ \ \ \ \ \ \ \ \ \ \ \ \left( \Pi _{j=1}^{n}\delta
_{(v_{j},x_{j},y_{j}:u_{w_{j}})}\right) \,\,E\left(
L_{w_{1}}^{u_{w_{1}}}...L_{w_{n}}^{u_{w_{n}}}\right) .$

$\square $
\end{proposition}

\strut

From now, rest of this paper, we will compute the $D_{G}$-valued cumulants
of the given $D_{G}$-valued random variable $a$. Let $w_{1},...,w_{n}\in
FP(G)$ be finite paths and \ $u_{1},...,u_{n}\in \{1,*\}$. Then, by the
M\"{o}bius inversion, we have

\strut

(2.2.1)$\ \ $

\begin{center}
$k_{n}\left( L_{w_{1}}^{u_{1}}~,...,~L_{w_{n}}^{u_{n}}\right) =\underset{\pi
\in NC(n)}{\sum }\widehat{E}(\pi )\left( L_{w_{1}}^{u_{1}}~\otimes
...\otimes ~L_{w_{n}}^{u_{n}}\right) \mu (\pi ,1_{n}),$
\end{center}

\strut

where $\widehat{E}=\left( E^{(n)}\right) _{n=1}^{\infty }$ is the moment
multiplicative bimodule map induced by the conditional expectation $E$ (See
[16]) and where $NC(n)$ is the collection of all noncrossing partition over $%
\{1,...,n\}.$ Notice that if $L_{w_{1}}^{u_{1}}...L_{w_{n}}^{u_{n}}$ does
not have the $*$-axis-property, then

\strut

\begin{center}
$E\left( L_{w_{1}}^{u_{1}}...L_{w_{n}}^{u_{n}}\right) =0_{D_{G}},$
\end{center}

\strut

by Section 2.1. Consider the noncrossing partition $\pi \in NC(n)$ with its
blocks $V_{1},...,V_{k}$. Choose one block $V_{j}=(j_{1},...,j_{k})\in \pi .$
Then we have that

\strut

(2.2.2) $\ $

\begin{center}
$\widehat{E}(\pi \mid _{V_{j}})\left( L_{w_{1}}^{u_{1}}~\otimes ...\otimes
~L_{w_{n}}^{u_{n}}\right) =E\left(
L_{w_{j_{1}}}^{u_{j_{1}}}d_{j_{2}}L_{w_{j_{2}}}^{u_{j_{2}}}...d_{j_{k}}L_{w_{j_{k}}}^{u_{j_{k}}}\right) , 
$
\end{center}

\strut where

\begin{center}
$d_{j_{i}}=\left\{ 
\begin{array}{ll}
1_{D_{G}} & 
\begin{array}{l}
\text{if there is no inner blocks} \\ 
\text{between }j_{i-1}\text{ \ and }j_{i}\text{ in }V_{j}
\end{array}
\\ 
&  \\ 
L_{v_{j_{i}}}\neq 1_{D_{G}} & 
\begin{array}{l}
\text{if there are inner blocks} \\ 
\text{between }j_{i-1}\text{ \ and }j_{i}\text{ in }V_{j},
\end{array}
\end{array}
\right. $
\end{center}

\strut

where $v_{j_{2}},...,v_{j_{k}}\in V(G).$ So, again by Section 2.1, $\widehat{%
E}(\pi \mid _{V_{j}})\left( L_{w_{1}}^{u_{1}}~\otimes ...\otimes
~L_{w_{n}}^{u_{n}}\right) $ is nonvanishing if and only if $%
L_{w_{j_{1}}}^{u_{j_{1}}}d_{j_{2}}L_{w_{j_{2}}}^{u_{j_{2}}}...d_{j_{k}}L_{w_{j_{k}}}^{u_{j_{k}}} 
$ has the $*$-axis-property, for all \ $j=1,...,n.$

\strut

Assume that

\begin{center}
$\widehat{E}(\pi \mid _{V_{j}})\left( L_{w_{1}}^{u_{1}}~\otimes ...\otimes
~L_{w_{n}}^{u_{n}}\right) =L_{v_{j}}$
\end{center}

and

\begin{center}
$\widehat{E}(\pi \mid _{V_{i}})\left( L_{w_{1}}^{u_{1}}~\otimes ...\otimes
~L_{w_{n}}^{u_{n}}\right) =L_{v_{i}}.$
\end{center}

\strut

If $v_{j}\neq v_{i},$ then the partition-dependent $D_{G}$-moment satisfies
that

\strut

\begin{center}
$\widehat{E}(\pi )\left( L_{w_{1}}^{u_{1}}~\otimes ...\otimes
~L_{w_{n}}^{u_{n}}\right) =0_{D_{G}}.$
\end{center}

\strut

This says that $\widehat{E}(\pi )\left( L_{w_{1}}^{u_{1}}~\otimes ...\otimes
~L_{w_{n}}^{u_{n}}\right) \neq 0_{D_{G}}$ if and only if there exists $v\in
V(G)$ such that

\strut

\begin{center}
$\widehat{E}(\pi \mid _{V_{j}})\left( L_{w_{1}}^{u_{1}}~\otimes ...\otimes
~L_{w_{n}}^{u_{n}}\right) =L_{v},$
\end{center}

\strut

for all \ $j=1,...,k.$

\strut

\begin{definition}
Let $NC(n)$ be the set of all \ noncrossing partition over $\{1,...,n\}$ and
let $L_{w_{1}}^{u_{1}},$ $...,$ $L_{w_{n}}^{u_{n}}\in \left(
W^{*}(G),E\right) $ be $D_{G}$-valued random variables, where $%
u_{1},...,u_{n}\in \{1,*\}.$ We say that the $D_{G}$-valued random variable $%
L_{w_{1}}^{u_{1}}...L_{w_{n}}^{u_{n}}$ is $\pi $-connected if the $\pi $%
-dependent $D_{G}$-moment of it is nonvanishing, for $\pi \in NC(n).$ In
other words, the random variable $L_{w_{1}}^{u_{1}}...L_{w_{n}}^{u_{n}}$ is $%
\pi $-connected, for $\pi \in NC(n),$ if

\strut

$\ \ \ \ \ \ \ \ \ \ \ \widehat{E}(\pi )\left( L_{w_{1}}^{u_{1}}~\otimes
...\otimes ~L_{w_{n}}^{u_{n}}\right) \neq 0_{D_{G}}.$

\strut

i.e., there exists a vertex $v\in V(G)$ such that

\strut

$\ \ \ \ \ \ \ \ \ \ \ \widehat{E}(\pi )\left( L_{w_{1}}^{u_{1}}~\otimes
...\otimes ~L_{w_{n}}^{u_{n}}\right) =L_{v}.$
\end{definition}

\strut

For convenience, we will define the following subset of $NC(n)$ ;

\strut

\begin{definition}
Let $NC(n)$ be the set of all noncrossing partitions over $\{1,...,n\}$ and
fix a $D_{G}$-valued random variable $L_{w_{1}}^{u_{1}}...L_{w_{n}}^{u_{n}}$
in $\left( W^{*}(G),E\right) ,$ where $u_{1},$ ..., $u_{n}\in \{1,*\}.$ For
the fixed $D_{G}$-valued random variable $%
L_{w_{1}}^{u_{1}}...L_{w_{n}}^{u_{n}},$define

\strut

$\ \ \ C_{w_{1},...,w_{n}}^{u_{1},...,u_{n}}\overset{def}{=}\{\pi \in
NC(n):L_{w_{1}}^{u_{1}}...L_{w_{n}}^{u_{n}}$ is $\pi $-connected$\},$

\strut

in $NC(n).$ Let $\mu $ be the M\"{o}bius function in the incidence algebra $%
I_{2}.$ Define the number $\mu _{w_{1},...,w_{n}}^{u_{1},...,u_{n}},$ for
the fixed $D_{G}$-valued random variable $%
L_{w_{1}}^{u_{1}}...L_{w_{n}}^{u_{n}},$ by

\strut

$\ \ \ \ \ \ \ \ \ \ \ \ \ \ \ \mu _{w_{1},...,w_{n}}^{u_{1},...,u_{n}}%
\overset{def}{=}\underset{\pi \in C_{w_{1},...,w_{n}}^{u_{1},...,u_{n}}}{%
\sum }\mu (\pi ,1_{n}).$
\end{definition}

\bigskip \strut

Assume that there exists $\pi \in NC(n)$ such that $%
L_{w_{1}}^{u_{1}}...L_{w_{n}}^{u_{n}}=L_{v}$ is $\pi $-connected. Then $\pi
\in C_{w_{1},...,w_{n}}^{u_{1},...,u_{n}}$ and there exists the maximal
partition $\pi _{0}\in C_{w_{1},...,w_{n}}^{u_{1},...,u_{n}}$ such that $%
L_{w_{1}}^{u_{1}}...L_{w_{n}}^{u_{n}}=L_{v}$ is $\pi _{0}$-connected.
(Recall that $NC(n)$ is a lattice. We can restrict this lattice ordering on $%
C_{w_{1},...,w_{n}}^{u_{1},...,u_{n}}$ and hence it is a POset, again.)
Notice that $1_{n}\in C_{w_{1},...,w_{n}}^{u_{1},...,u_{n}}.$ Therefore, the
maximal partition in $C_{w_{1},...,w_{n}}^{u_{1},...,u_{n}}$ is $1_{n}.$
Hence we have that ;

\bigskip \strut

\begin{lemma}
Let $L_{w_{1}}^{u_{1}}...L_{w_{n}}^{u_{n}}\in \left( W^{*}(G),E\right) $ be
a $D_{G}$-valued random variable having the $*$-axis-property. Then

\strut

$\ \ \ \ \ \ \ \ \ \ \ E\left( L_{w_{1}}^{u_{1}}...L_{w_{n}}^{u_{n}}~\right)
=\ \widehat{E}(\pi )\left( L_{w_{1}}^{u_{1}}~\otimes ...\otimes
~L_{w_{n}}^{u_{n}}\right) ,$

\strut

for all $\pi \in C_{w_{1},...,w_{n}}^{u_{1},...,u_{n}}.$
\end{lemma}

\bigskip

\begin{proof}
By the previous discussion, we can get the result.
\end{proof}

\strut \strut

By the previous lemmas, we have that

\strut

\begin{theorem}
Let $n\in 2\mathbb{N}$ and let $L_{w_{1}}^{u_{1}},...,L_{w_{n}}^{u_{n}}\in
\left( W^{*}(G),E\right) $ be $D_{G}$-valued random variables, where $%
w_{1},...,w_{n}\in FP(G)$ and $u_{j}\in \{1,*\},$ $j=1,...,n.$ Then

\strut

$\ \ \ \ \ \ \ k_{n}\left( L_{w_{1}}^{u_{1}}...L_{w_{n}}^{u_{n}}~\right)
=\mu _{w_{1},...,w_{n}}^{u_{1},...,u_{n}}\cdot
E(L_{w_{1}}^{u_{1}},...,L_{w_{n}}^{u_{n}}),$

\strut

where $\mu _{w_{1},...,w_{n}}^{u_{1},...,u_{n}}=\underset{\pi \in
C_{w_{1},...,w_{n}}^{u_{1},...,u_{n}}}{\sum }\mu (\pi ,1_{n}).$
\end{theorem}

\bigskip \strut \strut

\begin{proof}
We can compute that

\strut

$\ k_{n}\left( L_{w_{1}}^{u_{1}},\,...,L_{w_{n}}^{u_{n}}\right) =\underset{%
\pi \in NC(n)}{\sum }\widehat{E}(\pi )\left( L_{w_{1}}^{u_{1}}~\otimes
...\otimes ~L_{w_{n}}^{u_{n}}\right) \mu (\pi ,1_{n})$

\strut

$\ \ \ \ \ =\underset{\pi \in C_{w_{1},...,w_{n}}^{u_{1},...,u_{n}}}{\sum }%
\widehat{E}(\pi )\left( L_{w_{1}}^{u_{1}}~\otimes ...\otimes
~L_{w_{n}}^{u_{n}}\right) \mu (\pi ,1_{n})$

\strut

by the $\pi $-connectedness

\strut

$\ \ \ \ \ \ =\underset{\pi \in C_{w_{1},...,w_{n}}^{u_{1},...,u_{n}}}{\sum }%
E\left( L_{w_{1}}^{u_{1}}...~L_{w_{n}}^{u_{n}}\right) \mu (\pi ,1_{n})$

\strut

by the previous lemma

\strut

$\ \ \ \ \ \ =\left( \underset{\pi \in C_{w_{1},...,w_{n}}^{u_{1},...,u_{n}}%
}{\sum }\mu (\pi ,1_{n})\right) \cdot E\left(
L_{w_{1}}^{u_{1}}...~L_{w_{n}}^{u_{n}}\right) .$

\strut
\end{proof}

\strut

Now, we can get the following $D_{G}$-valued cumulants of the random
variable $a$ ;

\strut

\begin{corollary}
Let $n\in \mathbb{N}$ and let $a=a_{d}+a_{(*)}+a_{(non-*)}\in \left(
W^{*}(G),E\right) $ be our $D_{G}$-valued random variable. Then $k_{1}\left(
d_{1}a\right) =d_{1}a_{d}$ and $k_{n}\left( d_{1}a,...,d_{n}a\right)
=0_{D_{G}},$ for all odd $n.$ If $n\in \Bbb{N}\,\setminus \,\{1\},$ then

\strut \strut

$\ \ k_{n}\left( d_{1}a,...,d_{n}a\right) =~\underset{(v_{1},...,v_{n})\in
\Pi _{j=1}^{n}V(G:d_{j})}{\sum }\left( \Pi _{j=1}^{n}q_{v_{j}}\right) $

\strut

$\ \ \ \ \ \ \ \ \underset{(w_{1},...,w_{n})\in
FP_{*}(G:a)^{n},\,w_{j}=x_{j}w_{j}y_{j},\,u_{w_{j}}\in
\{1,*\},\,l_{w_{1},...,w_{n}}^{u_{1},...,u_{n}}\in LP_{n}^{*}}{\sum }\left(
\Pi _{j=1}^{n}p_{w_{j}}\right) $

\strut

$\ \ \ \ \ \ \ \ \ \ \left( \Pi _{j=1}^{n}\delta
_{(v_{j},x_{j},y_{j}:u_{j})}\right) \left( \mu
_{w_{1},...,w_{n}}^{u_{1},...,u_{n}}\cdot E\left(
L_{w_{1}}^{u_{w_{1}}}...L_{w_{n}}^{u_{w_{n}}}\right) \right) ,$

\strut

where $d_{1},...,d_{n}\in D_{G}$ are arbitrary. $\square $
\end{corollary}

\strut \strut \strut \strut

We have the following trivial $D_{G}$-valued moments and cumulants of an
arbitrary $D_{G}$-valued random variable ;

\strut

\begin{corollary}
Let $a\in \left( W^{*}(G),E\right) $ be a $D_{G}$-valued random variable and
let $n\in \mathbb{N}.$ Then

\strut

(1) The $n$-th trivial $D_{G}$-valued moment of $a$ is

\strut

$\ \ \ \ \ E(a^{n})=\,\underset{(w_{1},...,w_{n})\in \mathbb{F}%
^{+}(G:a)^{n},\,u_{w_{j}}\in
\{1,*\},\,l_{w_{1},...,w_{n}}^{u_{1},...,u_{n}}\in LP_{n}^{*}}{\sum }$

\strut

$\ \ \ \ \ \ \ \ \ \ \ \ \ \ \ \ \ \ \ \ \ \ \ \ \ \ \ \ \ \ \ \ \ \ \ \
\left( \Pi _{j=1}^{n}p_{w_{j}}\right) \cdot E\left(
L_{w_{1}}^{u_{1}},...,~L_{w_{n}}^{u_{n}}\right) .$

\strut

(2) The $n$-th trivial $D_{G}$-valued cumulant of $a$ is

\strut

$\ \ \ \ k_{1}\left( a\right) =E(a)=a_{d}$

\strut \strut

and

\strut

$\ \ \ \ \ k_{n}\left( \underset{n-times}{\underbrace{a,.....,a}}\right) =~%
\underset{(w_{1},...,w_{n})\in FP_{*}(G:a)^{n},\,\,\,u_{w_{j}}\in
\{1,*\},\,l_{w_{1},...,w_{n}}^{u_{1},...,u_{n}}\in LP_{n}^{*}}{\sum }$

\strut \strut

$\ \ \ \ \ \ \ \ \ \ \ \ \ \ \ \ \left( \Pi _{j=1}^{n}p_{w_{j}}\right)
\left( \mu _{w_{1},...,w_{n}}^{u_{1},...,u_{n}}\cdot E\left(
L_{w_{1}}^{u_{1}},...,~L_{w_{n}}^{u_{n}}\right) \right) ,$

\strut \strut

where $d_{1},...,d_{n}\in D_{G}$ are arbitrary. \ $\square $
\end{corollary}

\strut \bigskip

\strut

\strut

\section{\strut $D_{G}$-Freeness on $\left( W^{*}(G),E\right) $}

\strut

\strut

Like before, throughout this chapter, let $G$ be a countable directed graph
and $\left( W^{*}(G),E\right) $, the graph $W^{*}$-probability space over
its diagonal subalgebra $D_{G}.$ In this chapter, we will consider the $%
D_{G} $-valued freeness of given two random variables in $\left(
W^{*}(G),E\right) $. We will characterize the $D_{G}$-freeness of $D_{G}$%
-valued random variables $L_{w_{1}}$ and $L_{w_{2}},$ where $w_{1}\neq
w_{2}\in FP(G).$ And then we will observe the $D_{G}$-freeness of arbitrary
two $D_{G}$-valued random variables $a_{1}$ and $a_{2}$ in terms of their
supports. Let

\strut

(3.0) $\ a=\underset{w\in \mathbb{F}^{+}(G:a),\,u_{w}\in \{1,*\}}{\sum }%
p_{w}L_{w}^{u_{w}}$ \& $b=\underset{w^{\prime }\in \mathbb{F}%
^{+}(G:b),\,u_{w^{\prime }}\in \{1,*\}}{\sum }p_{w^{\prime }}L_{w^{\prime
}}^{u_{w^{\prime }}}$

\strut

be fixed $D_{G}$-valued random variables in $\left( W^{*}(G),E\right) $.

\strut

Now, fix $n\in \mathbb{N}$ and let $\left( a_{i_{1}}^{\varepsilon
_{i_{1}}},...,a_{i_{n}}^{\varepsilon _{i_{n}}}\right) \in
\{a,b,a^{*},b^{*}\}^{n},$ where $\varepsilon _{i_{j}}\in \{1,*\}.$ For
convenience, put

\strut

\begin{center}
$a_{i_{j}}^{\varepsilon _{i_{j}}}=\underset{w_{i_{j}}\in \mathbb{F}%
^{+}(G:a),\,\,u_{j}\in \{1,\ast \}}{\sum }p_{w_{j}}^{(j)}L_{w_{j}}^{u_{j}},$
for \ $j=1,...,n.$
\end{center}

\strut

Then, by the little modification of Section , we have that ;

\strut

(3.1)

\strut

$\ E\left( d_{i_{1}}a_{i_{1}}^{\varepsilon
_{i_{1}}}...d_{i_{n}}a_{i_{n}}^{\varepsilon _{i_{n}}}\right) $

\strut

$\ \ \ \ \ =\,\underset{(v_{i_{1}},...,v_{i_{n}})\in \Pi
_{k=1}^{n}V(G:d_{i_{k}})}{\sum }\left( \Pi _{k=1}^{n}q_{v_{i_{k}}}\right) $

\strut

$\ \ \ \ \ \ \ \ \ \ \ \ \ \ \ \ \ \underset{(w_{i_{1}},...,w_{i_{n}})\in
\Pi _{k=1}^{n}\mathbb{F}^{+}(G:a_{i_{k}}),%
\,w_{i_{j}}=x_{i_{j}}w_{i_{j}}y_{i_{j}},u_{i_{j}}\in \{1,\ast \}}{\sum }%
\left( \Pi _{k=1}^{n}p_{w_{i_{k}}}^{(k)}\right) $

\strut

$\ \ \ \ \ \ \ \ \ \ \ \ \ \ \ \ \ \ \ \ \ \ \ \ \ \ \ \ \ \ \ \ \ \ \left(
\Pi _{j=1}^{n}\delta _{(v_{i_{j}},x_{i_{j}},y_{i_{j}}:u_{i_{j}})}\right)
E\left( L_{w_{i_{1}}}^{u_{i_{1}}}...L_{w_{i_{n}}}^{u_{i_{n}}}\right) .$

\strut

\strut

Therefore, we have that

\strut

(3.2)

\strut

$\ \ k_{n}\left( d_{i_{1}}a_{i_{1}}^{\varepsilon
_{i_{1}}},...,d_{i_{n}}a_{i_{n}}^{\varepsilon _{i_{n}}}\right) =\underset{%
(v_{i_{1}},...,v_{i_{n}})\in \Pi _{k=1}^{n}V(G:d_{i_{k}})}{\sum }\left( \Pi
_{k=1}^{n}q_{v_{i_{k}}}\right) $

\strut

$\ \ \ \ \ \ \ \ \ \ \ \ \ \ \ \ \ \underset{(w_{i_{1}},...,w_{i_{n}})\in
\Pi
_{k=1}^{n}FP(G:a_{i_{k}}),\,w_{i_{j}}=x_{i_{j}}w_{i_{j}}y_{i_{j}},u_{i_{j}}%
\in \{1,\ast \}}{\sum }\left( \Pi _{k=1}^{n}p_{w_{i_{k}}}^{(k)}\right) $

\strut

$\ \ \ \ \ \ \ \ \ \ \ \ \ \ \ \ \ \ \ \ \left( \Pi _{j=1}^{n}\delta
_{(v_{i_{j}},x_{i_{j}},y_{i_{j}}:u_{i_{j}})}\right) \left( \mu
_{w_{i_{1}},...,w_{i_{n}}}^{u_{i_{1}},...,u_{i_{n}}}\cdot E\left(
L_{w_{i_{1}}}^{u_{i_{1}}}...L_{w_{i_{n}}}^{u_{i_{n}}}\right) \right) ,$

\strut

where $\mu _{w_{1},...,w_{n}}^{u_{1},...,u_{n}}=\underset{\pi \in
C_{w_{i_{1}},...,w_{i_{n}}}^{u_{i_{1}},...,u_{i_{n}}}}{\sum }\mu (\pi
,1_{n}) $ and

\bigskip

\begin{center}
$C_{w_{i_{1}},...,w_{i_{n}}}^{u_{i_{1}},...,u_{i_{n}}}=\{\pi \in
NC^{(even)}(n):L_{w_{1}}^{u_{1}}...L_{w_{n}}^{u_{n}}$ is $\pi $-connected$%
\}. $
\end{center}

\strut

So, we have the following proposition, by the straightforward computation ;

\bigskip \strut \strut

\begin{proposition}
Let $a,b\in \left( W^{*}(G),E\right) $ be $D_{G}$-valued random variables,
such that $a\notin W^{*}(\{b\},D_{G}),$ and let $\left(
a_{i_{1}}^{\varepsilon _{i_{1}}},...,a_{i_{n}}^{\varepsilon _{i_{n}}}\right)
\in \{a,b,a^{*},b^{*}\}^{n},$ for $n\in \mathbb{N}\setminus \{1\},$ where $%
\varepsilon _{i_{j}}\in \{1,*\},$ $j=1,...,n.$ Then\strut

\strut

\strut (3.3)

\strut

$\ \ \ \ \ k_{n}\left( d_{i_{1}}a_{i_{1}}^{\varepsilon
_{i_{1}}},...,d_{i_{n}}a_{i_{n}}^{\varepsilon _{i_{n}}}\right) $

\strut

$\ \ \ \ \ \ =\underset{(v_{1},...,v_{n})=(x,y,...,x,y)\in \Pi
_{j=1}^{n}V(G:d_{j})}{\sum }\left( \Pi _{j=1}^{n}q_{v_{j}}\right) $

\strut

$\ \ \ \ \ \ \ \ \ \ \underset{(w_{i_{1}},...,w_{i_{n}})\in \left( \Pi
_{k=1}^{n}FP_{*}(G:a_{i_{k}})\right) \cup
W_{*}^{iw_{1},...,w_{n}},\,w_{i_{j}}=x_{i_{j}}w_{i_{j}}y_{i_{j}}}{\sum }%
\left( \Pi _{k=1}^{n}p_{w_{i_{j}}}^{(k)}\right) $

\strut

$\ \ \ \ \ \ \ \ \ \ \ \ \ \ \left( \Pi _{j=1}^{n}\delta
_{(v_{i_{j}},x_{i_{j}},y_{i_{j}}:u_{i_{j}})}\right) \left( \mu
_{w_{i_{1}},...,w_{i_{n}}}^{u_{i_{1}},...,u_{i_{n}}}\cdot \Pr oj\left(
L_{w_{i_{1}}}^{u_{i_{1}}}...L_{w_{i_{n}}}^{u_{i_{n}}}\right) \right) $

\strut

\strut where $\mu _{n}=\underset{\pi \in
C_{w_{i_{1}},...,w_{i_{n}}}^{u_{i_{1}},...,u_{i_{n}}}}{\sum }\mu (\pi
,1_{n}) $ and

\strut

$\ \ \ \ \ \ W_{*}^{w_{1},...,w_{n}}=\{w\in FP_{*}^{c}(G:a)\cup
FP_{*}^{c}(G:b):$

$\ \ \ \ \ \ \ \ \ \ \ \ \ \ \ \ \ \ \ \ \ \ \ \ \ \ \ \ \ \ \ \ \ $both $%
L_{w}^{u_{w}}$ and $L_{w}^{u_{w}\,*}$\ are in $%
L_{w_{1}}^{u_{w_{1}}}...L_{w_{n}}^{u_{w_{n}}}\}.$

$\square $
\end{proposition}

\strut \strut \strut \strut \strut \strut

\begin{corollary}
Let $x$ and $y$ be the $D_{G}$-valued random variables in $\left(
W^{*}(G),E\right) $. The $D_{G}$-valued random variables $a$ and $b$ are
free over $D_{G}$ in $\left( W^{*}(G),E\right) $ if

\strut \strut \strut

$\ \ \ \ \ \ \ FP_{*}\left( G:P(x,x^{*})\right) \cap FP_{*}\left(
G:Q(y,y^{*})\right) =\emptyset $

and

$\ \ \ \ \ \ \ \ \ \ \ \ \ \ \ W_{*}^{\{P(x,x^{*}),\,Q(y,y^{*})\}}=\emptyset
,$

\strut

for all $P,Q\in \mathbb{C}[z_{1},z_{2}].$ \ $\square $
\end{corollary}

\strut \strut \strut \strut

By using (3.2), we can compute the mixed $D_{G}$-valued cumulants of two $%
D_{G}$-valued random variables. However, the formula is very abstract. So,
we will consider the above formula for fixed two generators of $W^{\ast }(G)$%
.

\strut

\begin{definition}
Let $G$ be a countable directed graph and $\mathbb{F}^{+}(G)$, the free
semigroupoid of $G$ and let $FP(G)$ be the subset of $\mathbb{F}^{+}(G)$
consisting of all finite paths. Define a subset $loop(G)$ of $FP(G)$
containing all loop finite paths or loops. (Remark that, in general, loop
finite paths are different from loop-edges. Clearly, all loop-edges are
loops in $FP(G).$) i.e.,

\strut

$\ \ \ \ \ \ \ loop(G)\overset{def}{=}\{l\in FP(G):l$ is a loop$\}\subset
FP(G).$

\strut

Also define the subset $loop^{c}(G)$ of $FP(G)$ consisting of all non-loop
finite path by

\strut

$\ \ \ \ \ \ \ \ \ loop^{c}(G)\overset{def}{=}FP(G)\,\setminus \,loop(G).$

\strut

Let $l\in loop(G)$ be a loop finite path. We say that $l$ is a \textbf{basic
loop} if there exists no loop $w\in loop(G)$ such that $l=w^{k},$ $k\in %
\mathbb{N}\,\setminus \,\{1\}.$ Define

\strut

$\ \ \ Loop(G)\overset{def}{=}\{l\in loop(G):l$ is a basic loop$%
\}\subsetneqq loop(G).$

\strut

Let $l_{1}=w_{1}^{k_{1}}$ and $l_{2}=w_{2}^{k_{2}}$ in $loop(G),$ where $%
w_{1},w_{2}\in Loop(G).$ We will say that the loops $l_{1}$ and $l_{2}$ are 
\textbf{diagram-distinct} if $w_{1}\neq w_{2}$ in $Loop(G).$ Otherwise, they
are not diagram-distinct.
\end{definition}

\strut

Now, we will introduce the more general diagram-distinctness of general
finite paths ;

\strut

\begin{definition}
(\textbf{Diagram-Distinctness}) We will say that the finite paths $w_{1}$
and $w_{2}$ are \textbf{diagram-distinct} if $w_{1}$ and $w_{2}$ have
different diagrams in the graph $G.$ Let $X_{1}$ and $X_{2}$ be subsets of $%
FP(G).$ The subsets $X_{1}$ and $X_{2}$ are said to be diagram-distinct if $%
x_{1}$ and $x_{2}$ are diagram-distinct for all pairs $(x_{1},x_{2})$ $\in $ 
$X_{1}\times X_{2}.$ This diagram-distinctness implies the
diagram-distinctness of loops.
\end{definition}

\strut \strut

Let $H$ be a directed graph with $V(H)=\{v_{1},v_{2}\}$ and $%
E(H)=\{e_{1}=v_{1}e_{1}v_{2},e_{2}=v_{2}e_{2}v_{1}\}.$ Then $l=e_{1}e_{2}$
is a loop in $FP(H)$ (i.e., $l\in loop(H)$). Moreover, it is a basic loop
(i.e., $l\in Loop(H)$). However, if we have a loop $%
w=e_{1}e_{2}e_{1}e_{2}=l^{2},$ then it is not a basic loop. i.e.,

\strut

\begin{center}
$l^{2}\in loop(H)\,\setminus \,Loop(H).$
\end{center}

\strut

If the graph $G$ contains at least one basic loop $l\in FP(G)$, then we have

\strut

\begin{center}
$\{l^{n}:n\in \mathbb{N}\}\subseteq loop(G)$ \ and \ $\{l\}\subseteq
Loop(G). $
\end{center}

\strut

Suppose that $l_{1}$ and $l_{2}$ are not diagram-distinct. Then, by
definition, there exists $w\in Loop(G)$ such that $l_{1}=w^{k_{1}}$ and $%
l_{2}=w^{k_{2}},$ for some $k_{1},k_{2}\in \mathbb{N}.$ On the graph $G,$
indeed, $l_{1}$ and $l_{2}$ make the same diagram. On the other hands, we
can see that if $w_{1}\neq w_{2}\in loop^{c}(G),$ then they are
automatically diagram-distinct.

\strut

\begin{lemma}
Suppose that $w_{1}\neq w_{2}\in loop^{c}(G)$ with $w_{1}=v_{11}w_{1}v_{12}$
and $w_{2}=v_{21}w_{2}v_{22}.$ Then $L_{w_{1}}$ and $L_{w_{2}}$ are free
over $D_{G}$ in $\left( W^{*}(G),E\right) .$
\end{lemma}

\strut

\begin{proof}
By definition, $L_{w_{1}}$ and $L_{w_{2}}$ are free over $D_{G}$ if and only
if all mixed $D_{G}$-valued cumulants of $W^{\ast }(\{L_{w_{1}}\},D_{G})$
and $W^{\ast }(\{L_{w_{2}}\},D_{G})$ vanish. Equivalently, all $D_{G}$%
-valued cumulants of $P\left( L_{w_{1}},L_{w_{1}}^{\ast }\right) $ and $%
Q\left( L_{w_{2}},L_{w_{2}}^{\ast }\right) $ vanish, for all $P,Q\in \mathbb{%
C}[z_{1},z_{2}].$ Since $w_{1}\neq w_{2}$ are non-loop edges, we can easily
verify that $w_{1}^{k_{1}}$ and $w_{2}^{k_{2}}$ are not admissible (i.e., $%
w_{1}^{k_{1}}\notin \mathbb{F}^{+}(G)$ and $w_{1}^{k_{2}}\notin \mathbb{F}%
^{+}(G)$), for all $k_{1},k_{2}\in \mathbb{N}\,\setminus \,\{1\}.$ This
shows that

\strut

$\ \ \ \ \ \ \ \ \ L_{w_{j}}^{k}=0_{D_{G}}=\left( L_{w_{j}}^{k\,}\right)
^{*},$ for $j=1,2.$

\strut

Thus, to show that $L_{w_{1}}$ and $L_{w_{2}}$ are free over $D_{G},$ it
suffices to show that all mixed $D_{G}$-valued cumulants of $P\left(
L_{w_{1}},L_{w_{1}}^{*}\right) $ and $Q\left( L_{w_{2}},L_{w_{2}}^{*}\right) 
$ vanish, for all $P,Q\in \mathbb{C}[z_{1},z_{2}]$ such that

\strut

$\ \ \ \ \ \ \ \ P(z_{1},z_{2})=\alpha _{1}z_{1}+\alpha
_{2}z_{1}z_{2}+\alpha _{3}z_{2}z_{1}+\alpha _{4}z_{2}$

and

$\ \ \ \ \ \ \ \ Q(z_{1},z_{2})=\beta _{1}z_{1}+\beta _{2}z_{1}z_{2}+\beta
_{3}z_{2}z_{1}+\beta _{4}z_{2},$

\strut \strut

where $\alpha ,\beta \in \mathbb{C}.$ So, for such $P$ and $Q,$ we have that

\strut

\ \ \ \ \ \ \strut $P\left( L_{w_{1}},L_{w_{1}}^{*}\right) =\alpha
_{1}L_{w_{1}}+\alpha _{2}L_{v_{11}}+\alpha _{3}L_{v_{12}}+\alpha
_{4}L_{w_{1}}^{*}$

and

$\ \ \ \ \ \ Q\left( L_{w_{2}},L_{w_{2}}^{*}\right) =\beta
_{1}L_{w_{2}}+\beta _{2}L_{v_{21}}+\beta _{3}L_{v_{22}}+\beta
_{4}L_{w_{2}}^{*}.$

\strut

Thus, we have that

\strut

$\ \ FP_{*}\left( G:P(L_{w_{1}},L_{w_{1}}^{*})\right) \supseteq \{w_{1}\},$
\ $FP_{*}\left( G:Q(L_{w_{2}},L_{w_{2}}^{*})\right) \supseteq \{w_{2}\}$

\strut

and

\strut

$\ \ FP_{*}^{c}\left( G:P(L_{w_{1}},L_{w_{1}}^{*})\right) \supseteq
\{w_{1}\},$ \ $FP_{*}^{c}\left( G:Q(L_{w_{2}},L_{w_{2}}^{*})\right)
\supseteq \{w_{2}\}.$

\strut

Remark that if $FP_{*}\left( G:P(L_{w_{1}},L_{w_{1}}^{*})\right) =\{w_{1}\},$
then $FP_{*}^{c}\left( G:P(L_{w_{1}},L_{w_{1}}^{*})\right) =\emptyset $, and
if$\ FP_{*}^{c}\left( G:P(L_{w_{1}},L_{w_{1}}^{*})\right) =\{w_{1}\},$ then $%
FP_{*}\left( G:P(L_{w_{1}},L_{w_{1}}^{*})\right) =\emptyset .$ The similar
relation holds for $Q\left( L_{w_{2}},L_{w_{2}}^{*}\right) .$ So, we have
that

\strut

$\ \ \ \ \ \ FP_{*}\left( G:P(L_{w_{1}},L_{w_{1}}^{*})\right) \cap
FP_{*}\left( G:Q(L_{w_{2}},L_{w_{2}}^{*})\right) =\emptyset $

and

$\ \ \ \ \ \ \ \ \ \ \ \ \
W_{*}^{\{P(L_{w_{1}},L_{w_{1}}^{*}),\,Q(L_{w_{2}},L_{w_{2}}^{*})\}}=%
\emptyset .$

\strut

Therefore, by the formula (3.4.3), we have the vanishing mixed $D_{G}$%
-valued cumulants of $P\left( L_{w_{1}},L_{w_{1}}^{\ast }\right) $ and $%
Q\left( L_{w_{2}},L_{w_{2}}^{\ast }\right) ,$ for all $n\in \mathbb{N}$ and
for all such $P,Q\in \mathbb{C}[z_{1},z_{2}].$ So, we can conclude that $%
L_{w_{1}}$ and $L_{w_{2}}$ are free over $D_{G}$ in $\left( W^{\ast
}(G),E\right) .$\strut
\end{proof}

\strut \strut \strut \strut \strut \strut \strut \strut \strut \strut

Now, we will consider the loop case.

\strut \strut

\begin{lemma}
Let $l_{1}\neq l_{2}\in Loop(G)$ be\textbf{\ basic loops} such that $%
l_{1}=v_{1}l_{1}v_{1}$ and $l_{2}=v_{2}l_{2}v_{2},$ for $v_{1},v_{2}\in V(G)$
(possibly $v_{1}=v_{2}$). i.e.e, two basic loops $l_{1}$ and $l_{2}$ are
diagram-distinct. Then the $D_{G}$-valued random variables $L_{l_{1}}$ and $%
L_{l_{2}}$ are free over $D_{G}$ in $\left( W^{*}(G),E\right) .$
\end{lemma}

\strut

\begin{proof}
Different from the non-loop case, if $l_{1}$ and $l_{2}$ are loops, then $%
l_{1}^{k_{1}}$ and $l_{2}^{k_{2}}$ exist in $FP(G),$ for all $k_{1},k_{2}\in %
\mathbb{N}.$ To show that $L_{l_{1}}$ and $L_{l_{2}}$ are free over $D_{G},$
it suffices to show that all mixed $D_{G}$-valued cumulants of $P\left(
L_{w_{1}},L_{w_{1}}^{\ast }\right) $ and $Q\left( L_{w_{2}},L_{w_{2}}^{\ast
}\right) $ vanish, for all $P,Q\in \mathbb{C}[z_{1},z_{2}].$ such that

\strut

$\ \ \ \ \ \ \ \ P(z_{1},z_{2})=f_{1}(z_{1})+f_{2}(z_{2})+P_{0}(z_{1},z_{2})$

and

$\ \ \ \ \ \ \ \
Q(z_{1},z_{2})=g_{1}(z_{1})+g_{2}(z_{2})+Q_{0}(z_{1},z_{2}), $

\strut \strut

where $f_{1},f_{2},g_{1},g_{2}\in \mathbb{C}[z]$ and $P_{0},Q_{0}\in \mathbb{%
C}[z_{1},z_{2}]$ such that $P_{0}$ and $Q_{0}$ does not contain polynomials
only in $z_{1}$ and $z_{2}.$ So, for such $P$ and $Q,$ we have that

\strut

\ \ \ \ \ \ \ \strut $P\left( L_{l_{1}},L_{l_{1}}^{*}\right)
=f_{1}(L_{l_{1}})+f_{2}(L_{l_{1}}^{*})+P_{0}(L_{l_{1}},L_{l_{1}}^{*})$

and

$\ \ \ \ \ \ \ Q\left( L_{l_{2}},L_{l_{2}}^{*}\right)
=g_{1}(L_{l_{2}})+g_{2}(L_{l_{2}}^{*})+Q_{0}(L_{l_{2}},L_{l_{2}}^{*}).$

\strut

Notice that $L_{l_{j}}^{k}=L_{l_{j}^{k}},$ for all $k\in \mathbb{N},$ \ $%
j=1,2.$ Also, notice that

\strut

$\ \ \ \ \ \ \ \ \ P_{0}\left( L_{l_{1}},L_{l_{1}}^{*}\right)
=f_{1}^{0}(L_{l_{1}})+f_{2}^{0}(L_{l_{1}}^{*})+\alpha L_{v_{1}}$

and

$\ \ \ \ \ \ \ \ \ Q_{0}\left( L_{l_{2}},L_{l_{2}}^{*}\right)
=g_{1}^{0}(L_{l_{2}})+g_{2}^{0}(L_{l_{2}}^{*})+\beta L_{v_{2}},$

\strut

where $f_{1}^{0},f_{2}^{0},g_{1}^{0},g_{2}^{0}\in \mathbb{C}[z]$ and $\alpha
,\beta \in \mathbb{C}$, by the fact that

\strut

\ \ \ \ \ \ \ \ \ \ \ \ \ \ \ \ \ $\ \
L_{l_{j}}^{*}L_{l_{j}}=L_{v_{j}}=L_{l_{j}}L_{l_{j}}^{*},$

\strut

under the weak-topology. So, finally, we have that

\strut

$\ \ \ \ \ \ \ \ \ P\left( L_{l_{1}},L_{l_{1}}^{*}\right) =\mathbf{f}%
_{1}(L_{l_{1}})+\mathbf{f}_{2}(L_{l_{1}}^{*})+\alpha L_{v_{1}}$

and

$\ \ \ \ \ \ \ \ \ Q\left( L_{l_{2}},L_{l_{2}}^{*}\right) =\mathbf{g}%
_{1}(L_{l_{2}})+\mathbf{g}_{1}(L_{l_{2}}^{*})+\beta L_{v_{2}},$

\strut

where $\mathbf{f}_{1},\mathbf{f}_{2},\mathbf{g}_{1},\mathbf{g}_{2}\in %
\mathbb{C}[z]$ and $\alpha ,\beta \in \mathbb{C}.$ Thus, we have that

\strut

$\ \ FP_{*}\left( G:P(L_{l_{1}},L_{l_{1}}^{*})\right) \subseteq
\{l_{1}^{k}\}_{k=1}^{\infty },$ \ $FP_{*}\left(
G:Q(L_{l_{2}},L_{l_{2}}^{*})\right) \subseteq \{l_{2}^{k}\}_{k=1}^{\infty }$

\strut

and

\strut

$\ \ FP_{*}^{c}\left( G:P(L_{l_{1}},L_{l_{1}}^{*})\right) \subseteq
\{l_{1}^{k}\}_{k=1}^{\infty },$ \ $FP_{*}^{c}\left(
G:Q(L_{l_{2}},L_{l_{2}}^{*})\right) \subseteq \{l_{2}^{k}\}_{k=1}^{\infty }.$

\strut

So, we have that

\strut

$\ \ \ \ \ \ FP_{*}\left( G:P(L_{w_{1}},L_{w_{1}}^{*})\right) \cap
FP_{*}\left( G:Q(L_{w_{2}},L_{w_{2}}^{*})\right) =\emptyset ,$

\strut

because $l_{1}$ and $l_{2}$ are in $Loop(G)$ (and hence if $l_{1}\neq l_{2},$
then they are diagram-distinct.) And we have that

\strut

$\ \ \ \ \ \ \ \ \ \ \ \ \
W_{*}^{\{P(L_{w_{1}},L_{w_{1}}^{*}),\,Q(L_{w_{2}},L_{w_{2}}^{*})\}}=%
\emptyset .$

\strut

Therefore, by the formula (3.4.3), we have the vanishing mixed $D_{G}$%
-valued cumulants of $P\left( L_{l_{1}},L_{l_{1}}^{\ast }\right) $ and $%
Q\left( L_{l_{2}},L_{l_{2}}^{\ast }\right) ,$ for all $n\in \mathbb{N}$ and
for all $P,Q\in \mathbb{C}[z_{1},z_{2}].$ Since $P$ and $Q$ are arbitrary,
we can conclude that $L_{l_{1}}$ and $L_{l_{2}}$ are free over $D_{G}$ in $%
\left( W^{\ast }(G),E\right) .$\strut
\end{proof}

\strut

Notice that we assumed that the loops $l_{1}$ and $l_{2}$ are basic loops in
the previous lemma. Since they are distinct basic loops, they are
automatically diagram-distinct. Now, assume that $l_{1}$ and $l_{2}$ are not
diagram-distinct. i.e., there exists a basic loop $w\in Loop(G)$ such that $%
l_{1}=w^{k_{1}}$ and $l_{2}=w^{k_{2}},$ for some $k_{1},k_{2}\in \mathbb{N}.$
In other words, the loops $l_{1}$ and $l_{2}$ have the same diagram in the
graph $G.$ Then the $D_{G}$-valued random variables $L_{l_{1}}$ and $%
L_{l_{2}}$ are not free over $D_{G}$ in $\left( W^{\ast }(G),E\right) $. See
the next example ;

\strut

\begin{example}
Let $G_{1}$ be a directed graph with $V(G_{1})=\{v\}$ and $%
E(G_{1})=\{l=vlv\}.$ So, in this case,

\strut

$\ \ \ \ \ \ \ E(G_{1})=Loop(G_{1}),$ \ $FP(G_{1})=loop(G_{1}),$

\strut and

$\ \ \ \ \ \ \ \ \ \ \ loop(G_{1})=\{l^{k}:k\in \mathbb{N}\}.$

\strut

Thus, even if $w_{1}\neq w_{2}\in loop(G_{1}),$ $w_{1}$ and $w_{2}$ are Not
diagram-distinct. Take $l^{2}$ and $l^{3}$ in $FP(G_{1}).$ Then the $%
D_{G_{1}}$-valued random variable $L_{l^{2}}$ and $L_{l^{3}}$ are not free
over $D_{G_{1}}$ in $\left( W^{*}(G_{1}),E\right) .$ Indeed, let's take $%
P,Q\in \mathbb{C}[z_{1},z_{2}]$ as

\strut

$\ \ \ \ \ \ \ P(z_{1},z_{2})=z_{1}^{3}+z_{2}^{3}$ \ \ and \ \ $%
Q(z_{1},z_{2})=z_{1}^{2}+z_{2}^{2}.$

\strut Then

$\ \ \ \ \ \ \ \ \ P\left( L_{l^{2}},L_{l^{2}}^{*}\right)
=L_{l^{2}}^{3}+L_{l^{2}}^{*\,3}=L_{l^{6}}+L_{l^{6}}^{*}$

and

\ \ \ \ \ \ \ \ \ $Q\left( L_{l^{3}},L_{l^{3}}^{*}\right)
=L_{l^{3}}^{2}+L_{l^{3}}^{*\,\,2}=L_{l^{6}}+L_{l^{6}}^{*}.$

\strut

Then

\strut

$\ k_{2}\left(
P(L_{l^{2}},L_{l^{2}}^{*}),\,Q(L_{l^{3}},L_{l^{3}}^{*})\right) =k_{2}\left(
L_{l^{6}}+L_{l^{6}}^{*},\,L_{l^{6}}+L_{l^{6}}^{*}\right) $

\strut

$\ \ \ =\mu _{l^{6},l^{6}}^{1,*}\Pr oj\left( L_{l^{6}},L_{l^{6}}^{*}\right)
+\mu _{l^{6},l^{6}}^{*,1}\Pr oj\left( L_{l^{6}}^{*},L_{l^{6}}\right) $

\strut

$\ \ \ =\mu _{l^{6},l^{6}}^{1,*}L_{v}+\mu _{l^{6},l^{6}}^{*,1}L_{v}=\left(
\mu _{l^{6},l^{6}}^{1,*}+\mu _{l^{6},l^{6}}^{*,1}\right) L_{v}$

\strut

$\ \ \ =2L_{v}\neq 0_{D_{G}},$

\strut

since $\mu _{l^{6},l^{6}}^{1,*}=\mu (1_{2},1_{2})=1=\mu
_{l^{6},l^{6}}^{*,1}. $ This says that there exists at least one
nonvanishing mixed $D_{G}$-valued cumulant of $W^{*}\left(
\{L_{l^{3}}\},D_{G_{1}}\right) $ and $W^{*}\left(
\{L_{l^{2}}\},D_{G_{1}}\right) .$ Therefore, $L_{l^{3}}$ and $L_{l^{2}}$ are
not free over $D_{G_{1}}$ in $\left( W^{*}(G_{1}),E\right) .$
\end{example}

\strut

As we have seen before, if two loops $l_{1}$ and $l_{2}$ are not
diagram-distinct, then $D_{G}$-valued random variables $L_{l_{1}}$ and $%
L_{l_{2}}$ are Not free over $D_{G}.$ However, if $l_{1}$ and $l_{2}$ are
diagram-distinct, we can have the following lemma, by the previous lemma ;

\strut

\begin{lemma}
Let $l_{1}\neq l_{2}\in loop(G)$ be loops and assume that $%
l_{1}=w_{1}^{k_{1}}$ and $l_{2}=w_{2}^{k_{2}},$ where $w_{1},w_{2}\in
Loop(G) $ are basic loops and $k_{1},k_{2}\in \mathbb{N}.$ If $w_{1}\neq
w_{2}\in Loop(G),$ then the $D_{G}$-valued random variables $L_{l_{1}}$ and $%
L_{l_{2}} $ are free over $D_{G}$ in $\left( W^{*}(G),E\right) .$ $\square $
\end{lemma}

\strut \strut \strut \strut

Finally, we will observe the following case when we have a loop and a
non-loop finite path.

\strut

\begin{lemma}
Let $l\in loop(G)$ and $w\in loop^{c}(G).$ Then the $D_{G}$-valued random
variables $L_{l}$ and $L_{w}$ are free over $D_{G}$ in $\left(
W^{*}(G),E\right) .$
\end{lemma}

\strut

\begin{proof}
let $l\in loop(G)$ and $w\in loop^{c}(G)$ and let $L_{l}$ and $L_{w}$ be the
corresponding $D_{G}$-valued random variables in $\left( W^{\ast
}(G),E\right) .$ Then, for all $P,Q\in \mathbb{C}[z_{1},z_{2}],$ we have that

\strut

$\ \ \ \ \ \ \ \ FP_{*}\left( G:P(L_{l},L_{l}^{*})\right) \cap FP_{*}\left(
G:Q(L_{w},L_{w}^{*})\right) =\emptyset ,$

\strut

\strut since

$\ \ \ \ \ \ \ \ \ FP_{*}\left( G:P(L_{l},L_{l}^{*})\right) \subseteq
\{l^{k}:k\in \mathbb{N}\}\subset loop(G)$

and

$\ \ \ \ \ \ \ \ \ \ \ \ \ FP_{*}\left( G:Q(L_{w},L_{w}^{*})\right)
=\{w\}\subset loop^{c}(G).$

\strut

Also, since $loop(G)\cap loop^{c}(G)=\emptyset ,$ we can get that

\strut

$\ \ \ \ \ \ \ \ \ \ \ \ \ \ \ \ \
W_{*}^{\{P(L_{l},L_{l}^{*}),\,Q(L_{w},L_{w}^{*})\}}=\emptyset ,$

\strut

for all $P,Q\in \mathbb{C}[z_{1},z_{2}].$ Therefore, the $D_{G}$-valued
random variables $L_{l}$ and $L_{w}$ are free over $D_{G}$ in $\left(
W^{*}(G),E\right) .$
\end{proof}

\strut

\strut Now, we can summarize the above lemmas in this section as follows and
this theorem is one of the main result of this paper. The theorem is the
characterization of $D_{G}$-freeness of generators of $W^{*}(G)$ over $D_{G}$%
.

\strut \strut

\begin{theorem}
Let $w_{1},w_{2}\in FP(G)$ be finite paths. The $D_{G}$-valued random
variables $L_{w_{1}}$ and $L_{w_{2}}$ in $\left( W^{*}(G),E\right) $ are
free over $D_{G}$ if and only if $w_{1}$ and $w_{2}$ are diagram-distinct.
\end{theorem}

\strut

\begin{proof}
$\Leftarrow $) Suppose that finite paths $w_{1}$ and $w_{2}$ are
diagram-distinct. Then the $D_{G}$-valued random variables $L_{w_{1}}$ and $%
L_{w_{2}}$ are free over $D_{G},$ by the previous lemmas.

\strut

($\Rightarrow $) Let $L_{w_{1}}$ and $L_{w_{2}}$ be free over $D_{G}$ in $%
\left( W^{*}(G),E\right) .$ Now, assume that $w_{1}$ and $w_{2}$ are not
diagram-distinct. We will observe the following cases ;

\strut

(Case I) The finite paths $w_{1},w_{2}\in loop(G).$ Since they are not
diagram-distinct, there exists a basic loop $l$ $=$ $vlv$ $\in $ $Loop(G),$
with $v\in V(G),$ such that $w_{1}=l^{k_{1}}$ and $w_{2}=l^{k_{2}},$ for
some $k_{1},k_{2}\in \mathbb{N}.$ As we have seen before, $L_{w_{1}}$ and $%
L_{w_{2}}$ are not free over $D_{G}$ in $\left( W^{*}(G),E\right) .$ Indeed,
if we let $k\in \mathbb{N}$ such that $k_{1}\mid k$ \ and \ $k_{2}\mid k$ \
with $k=k_{1}n_{1}=k_{2}n_{2},$ for $n_{1},n_{2}\in \mathbb{N},$ then we can
take $P,Q\in \mathbb{C}[z_{1},z_{2}]$ defined by

\strut

$\ \ \ \ \ \ P(z_{1},z_{2})=z_{1}^{n_{1}}+z_{2}^{n_{1}}$ \ and \ $%
Q(z_{1},z_{2})=z_{1}^{n_{2}}+z_{2}^{n_{2}}.$

\strut

And then

\strut

$\ \ \ \ P\left( L_{w_{1}},L_{w_{1}}^{*}\right)
=L_{w_{1}}^{n_{1}}+L_{w_{1}}^{*\,%
\,n_{1}}=L_{l^{k_{1}}}^{n_{1}}+L_{l^{k_{1}}}^{*\,%
\,n_{1}}=L_{l^{k}}+L_{l^{k}}^{*}$

\strut

and

\strut

$\ \ \ \ Q\left( L_{w_{2}},L_{w_{2}}^{*}\right)
=L_{w_{2}}^{n_{2}}+L_{w_{2}}^{*\,%
\,n_{2}}=L_{l^{k_{2}}}^{n_{2}}+L_{l^{k_{2}}}^{*%
\,n_{2}}=L_{l^{k}}+L_{l^{k}}^{*}.$

\strut

So,

\strut

$\ \ \ 
\begin{array}{ll}
k_{2}\left( P(L_{w_{1}},L_{w_{1}}^{*}),\,Q(L_{w_{2}},L_{w_{2}}^{*})\right) & 
=k_{2}\left( L_{l^{k}}+L_{l^{k}}^{*},\,L_{l^{k}}+L_{l^{k}}^{*}\right) \\ 
&  \\ 
& =\mu _{l^{k},l^{k}}^{1,*}L_{v}+\mu _{l^{k},l^{k}}^{*,1}L_{v} \\ 
&  \\ 
& =2L_{v}\neq 0_{D_{G}}.
\end{array}
$

\strut

Therefore, $P\left( L_{w_{1}},L_{w_{1}}^{*}\right) $ and $Q\left(
L_{w_{2}},L_{w_{2}}^{*}\right) $ are not free over $D_{G}.$ This shows that $%
W^{*}\left( \{L_{w_{1}}\},D_{G}\right) $ and $W^{*}\left(
\{L_{w_{2}}\},D_{G}\right) $ are not free over $D_{G}$ in $\left(
W^{*}(G),E\right) $ and hence $L_{w_{1}}$ and $L_{w_{2}}$ are not free over $%
D_{G}.$ This contradict our assumption.

\strut

(Case II) Suppose that the finite paths $w_{1},w_{2}$ are non-loop finite
paths in $loop^{c}(G)$ and assume that they are not diagram-distinct. Since
they are not diagram-distinct, they are identically equal. Therefore, they
are not free over $D_{G}$ in $\left( W^{*}(G),E\right) .$

\strut

(Case III) Let $w_{1}\in loop(G)$ and $w_{2}\in loop^{c}(G).$ They are
always diagram-distinct.

\strut

Let $L_{w_{1}}$ and $L_{w_{2}}$ are free over $D_{G}$ and assume that $w_{1}$
and $w_{2}$ are not diagram-distinct. Then $L_{w_{1}}$ and $L_{w_{2}}$ are
not free over $D_{G}$, by the Case I, II and III. So, this contradict our
assumption.
\end{proof}

\strut

The previous theorem characterize the $D_{G}$-freeness of two partial
isometries $L_{w_{1}}$ and $L_{w_{2}},$ where $w_{1},w_{2}\in FP(G).$ This
characterization shows us that the diagram-distinctness of finite paths
determine the $D_{G}$-freeness of corresponding creation operators.

\strut

Let $a$ and $b$ be the given $D_{G}$-valued random variables in (3.0). We
can get the necessary condition for the $D_{G}$-freeness of $a$ and $b,$ in
terms of their supports. \strut \strut \strut Recall that we say that the
two subsets $X_{1}$ and $X_{2}$ of $FP(G)$ are said to be diagram-distinct
if $x_{1}$ and $x_{2}$ are diagram-distinct, for all pairs $(x_{1},x_{2})$ $%
\in $ $X_{1}$ $\times $ $X_{2}.$

\strut \strut

\begin{theorem}
Let $a,b\in \left( W^{*}(G),E\right) $ be $D_{G}$-valued random variables
with their supports $\mathbb{F}^{+}(G:a)$ and $\mathbb{F}^{+}(G:b).$ The $%
D_{G}$-valued random variables $a$ and $b$ are free over $D_{G}$ in $\left(
W^{*}(G),E\right) $ if $FP(G:a_{1})$ and $FP(G:a_{2})$ are diagram-distinct.
\end{theorem}

\strut

\begin{proof}
For convenience, let's denote $a$ and $b$ by $a_{1}$ and $a_{2},$
respectively. Assume that the supports of $a_{1}$ and $a_{2},$ $\mathbb{F}%
^{+}(G:a_{1})$ and $\mathbb{F}^{+}(G:a_{2})$ are diagram-distinct. Then by
the previous $D_{G}$-freeness characterization,

\strut

$\ \ \ \underset{l\in FP(G:a_{1}),\,u_{l}\in \{1,*\}}{\sum }%
p_{l}^{(1)}L_{l}^{u_{l}}$ \ and \ $\underset{l\in FP(G:a_{2}),\,u_{l}\in
\{1,*\}}{\sum }p_{l}^{(2)}L_{l}^{u_{l}}$

\strut

are free over $D_{G}$ in $\left( W^{*}(G),E\right) .$ Indeed, since $%
FP(G:a_{1})$ and $FP(G:a_{2})$ are diagram-distinct, all summands $L_{w_{1}}$%
's of $a_{1}$ and $L_{w_{2}}$'s of $a_{2}$ are free over $D_{G}$ in $\left(
W^{*}(G),E\right) .$ Therefore, $a_{1}$ and $a_{2}$ are free over $D_{G}$ in 
$\left( W^{*}(G),E\right) .$\strut \strut \strut
\end{proof}

\strut \strut \strut

\strut \strut \strut

\strut \strut \strut

\section{Examples}

\strut \strut

\strut

\strut

In this chapter, as examples, we will compute the trivial $D_{G}$-valued
moments and cumulants of the generating operator $T$ of the graph $W^{*}$%
-algebra $W^{*}(G).$ Let $G$ be a countable directed graph and let $\left(
W^{*}(G),E\right) $ be the graph $W^{*}$-probability space over its diagonal
subalgebra. Let $a\in \left( W^{*}(G),E\right) $ be a $D_{G}$-valued random
variable. Recall that the trivial $D_{G}$-valued $n$-th moments and
cumulants of $a$ are defined by

\strut

\begin{center}
$E(a^{n})$ \ \ and \ $k_{n}\left( \underset{n\text{-times}}{\underbrace{%
a,.......,a}}\right) .$
\end{center}

\strut

In this chapter, we will deal with the following special $D_{G}$-valued
random variable ;

\strut

\begin{definition}
Define an operator $T$ in $W^{*}(G)$ by

\strut

(5.1)$\strut $ $\ \ \ \ \ \ \ \ \ \ \ \ \ T=\underset{e\in E(G)}{\sum }%
\left( L_{e}+L_{e}^{*}\right) .$

\strut

We will call $T$ the generating operator of $W^{*}(G).$ The self-adjoint
operators $L_{e}+L_{e}^{*},$ for $e\in E(G),$ are called the block operators
of $T.$
\end{definition}

\strut

\begin{example}
Let $G$ be a one-vertex directed graph with $N$-edges. i.e.,

\strut

$\ \ \ \ V(G)=\{v\}$ \ and \ $E(G)=\{e_{j}=ve_{j}v:j=1,...,N\}.$

\strut

Then the graph $W^{*}$\strut -algebra $W^{*}(G)$ satisfies that

\strut

(5.2)$\ \ \ \ W^{*}(G)=D_{G}*_{D_{G}}\left( \underset{j=1}{\overset{N}{%
\,*_{D_{G}}}}\left( W^{*}(\{L_{e_{j}}\},D_{G})\right) \right) ,$

\strut

by Chapter 4. Notice that $D_{G}=\Delta _{1}=\mathbb{C}.$ Therefore, the
formula (5.2) is rewritten by

\strut

(5.4) \ \ \ $\ \ \ \ \ \ \ W^{*}(G)=\,\underset{j=1}{\overset{N}{*}}\left(
W^{*}(\{L_{e_{j}}\})\right) ,$

\strut

where $*$ means the usual (scalar-valued) free product of Voiculescu. Also
notice that $1_{D_{G}}=L_{v}=1\in \mathbb{C}$ and

\strut

(5.5) $\ \ \ L_{e_{j}}^{*}L_{e_{j}}=L_{v}=1=L_{e_{j}}L_{e_{j}}^{*},$ for all 
$j=1,...,N.$

\strut

This shows that $L_{e_{j}}$'s are unitary in $W^{*}(G),$ for all $j=1,...,N.$
Now, define the generating operator $T=\sum_{j=1}^{N}\left(
L_{e_{j}}+L_{e_{j}}^{*}\right) $ of $W^{*}(G).$ It is easy to see that each
block operator $L_{e_{j}}+L_{e_{j}}^{*}$ is semicircular, by Voiculescu, for
all $j=1,...,N.$ (Remember the construction of creation operators $L_{e_{j}}$%
's and see [9].) Futhermore, by Chapter 3, we can get that all blocks $%
L_{e_{j}}+L_{e_{j}}^{*}$'s are free from each other in the graph $W^{*}$%
-probability space $\left( W^{*}(G),E\right) .$

\strut

By (5.4), the canonical conditional expectation $E:W^{*}(G)\rightarrow D_{G}$
is the faithful linear functional. Moreover, by (5.5), this linear
functional $E$ is a trace in the sense that $E(ab)=E(ba),$ for all $a,b\in
W^{*}(G).$ From now, to emphasize that $E$ is a trace, we will denote $E$ by 
$tr.$

\strut

Let's compute the $n$-th cumulant of $T$ ;

\strut

\strut (5.6)

$\ 
\begin{array}{ll}
k_{n}\left( T,...,T\right) & =k_{n}\left( \sum_{j=1}^{N}\left(
L_{e_{j}}+L_{e_{j}}^{*}\right) ,...,\sum_{j=1}^{N}\left(
L_{e_{j}}+L_{e_{j}}^{*}\right) \right) \\ 
&  \\ 
& =\sum_{j=1}^{N}k_{n}\left(
(L_{e_{j}}+L_{e_{j}}^{*}),...,(L_{e_{j}}+L_{e_{j}}^{*})\right) ,
\end{array}
$ $\ $

\strut \strut

by the mutual freeness of $\{L_{e_{j}},L_{e_{j}}^{*}\}$'s on $\left(
W^{*}(G),tr\right) $, for $j=1,...,N.$ Observe that

\strut

$\ \ \ k_{n}\left(
(L_{e_{j}}+L_{e_{j}}^{*}),...,(L_{e_{j}}+L_{e_{j}}^{*})\right) $

\strut

(5.7) $\ \ \ =\left\{ 
\begin{array}{ll}
k_{2}\left( (L_{e_{j}}+L_{e_{j}}^{*}),(L_{e_{j}}+L_{e_{j}}^{*})\right) & 
\text{if }n=2 \\ 
&  \\ 
0 & \text{otherwise,}
\end{array}
\right. $

\strut

by the semicircularity of $L_{e_{j}}+L_{e_{j}}^{*},$ for $j=1,...,N.$ By
(5.7), the formula (5.6) is

\strut

$\ \ \ k_{n}\left( T,...,T\right) $

\strut

(5.8) $\ \ \ \ =\left\{ 
\begin{array}{ll}
\sum_{j=1}^{N}k_{2}\left(
L_{e_{j}}+L_{e_{j}}^{*},L_{e_{j}}+L_{e_{j}}^{*}\right) & \text{if }n=2 \\ 
&  \\ 
0 & \text{otherwise}
\end{array}
\right. $

$\strut \ \ \ \ $

\strut Now, observe $k_{2}\left(
L_{e_{j}}+L_{e_{j}}^{*},L_{e_{j}}+L_{e_{j}}^{*}\right) $ ;

\strut

$\ k_{2}\left( L_{e_{j}}+L_{e_{j}}^{*},L_{e_{j}}+L_{e_{j}}^{*}\right) $

\strut

$\ \ \ =k_{2}\left( L_{e_{j}},L_{e_{j}}\right) +k_{2}\left(
L_{e_{j}},L_{e_{j}}^{*}\right) +k_{2}\left( L_{e_{j}}^{*},L_{e_{j}}\right)
+k_{2}\left( L_{e_{j}}^{*},L_{e_{j}}^{*}\right) $

\strut

$\ \ \ =0+k_{2}\left( L_{e_{j}},L_{e_{j}}^{*}\right) +k_{2}\left(
L_{e_{j}}^{*},L_{e_{j}}\right) +0$

\strut

by Section 2.1

\strut

$\ \ \ =tr\left( L_{e_{j}}L_{e_{j}}^{*}\right) +tr\left(
L_{e_{j}}^{*}L_{e_{j}}\right) =2\cdot tr\left( L_{e_{j}}^{*}L_{e_{j}}\right) 
$

\strut

since $tr$ is a trace

\strut

$\ \ \ =2\cdot L_{v}=2,$

\strut

for $j=1,...,N,$ by Section 2.1 and 2.2. So, we can get that

\strut

(5.9) \ $k_{n}\left( T,...,T\right) =\left\{ 
\begin{array}{lll}
2N &  & \text{if }n=2 \\ 
&  &  \\ 
0 &  & \text{otherwise.}
\end{array}
\right. $

\strut

Now, we can compute the trivial moments of $T,$ via the M\"{o}bius inversion.

\strut

$\ tr\left( T^{n}\right) =\underset{\pi \in NC(n)}{\sum }k_{\pi }\left(
a,...,a\right) $

\strut

where $k_{\pi }(a,...,a)=\underset{V\in \pi }{\prod }k_{\left| V\right|
}\left( \underset{\left| V\right| \text{-times}}{\underbrace{a,.......,a}}%
\right) ,$ for each $\pi \in NC(n),$ by Nica and Speicher (See [1] and [17])

\strut

$\ \ \ \ \ \ =\underset{\pi \in NC_{2}(n)}{\sum }k_{\pi }(a,...,a)=\underset{%
\pi \in NC_{2}(n)}{\sum }\,\left( \underset{V\in \pi }{\prod }k_{\left|
V\right| }\left( a,a\right) \right) $

\strut

where $NC_{2}(n)=\{\pi \in NC(n):V\in \pi \Leftrightarrow \left| V\right|
=2\}$ is the collection of all noncrossing pairings

\strut

(5.10) $\ $

$\ \ \ \ \ \ =\underset{\pi \in NC_{2}(n)}{\sum }\,\left( \underset{V\in \pi 
}{\prod }2N\right) =\underset{\pi \in NC_{2}(n)}{\sum }\left( 2N\right)
^{\left| \pi \right| },$

\strut

where $\left| \pi \right| \overset{def}{=}$ the number of blocks in $\pi .$
Notice that the above formula (5.10) shows us that the $n$ should be even,
because $NC_{2}(n)$ is nonempty when $n$ is even. Therefor,

\strut

(5.11) \ \ \ $tr\left( T^{n}\right) =\left\{ 
\begin{array}{ll}
\underset{\pi \in NC_{2}(n)}{\sum }\left( 2N\right) ^{\left| \pi \right| } & 
\text{if }n\text{ is even} \\ 
&  \\ 
0 & \text{if }n\text{ is odd.}
\end{array}
\right. $

\strut

Also, notice that if $\pi \in NC_{2}(n),$ then $\left| \pi \right| =\frac{n}{%
2},$ for all even number $n\in \mathbb{N}.$ So,

\strut

\ \ \ $tr\left( T^{n}\right) =\left\{ 
\begin{array}{ll}
\left| NC_{2}(n)\right| \cdot \left( 2N\right) ^{\frac{n}{2}} & \text{if }n%
\text{ is even} \\ 
&  \\ 
0 & \text{if }n\text{ is odd}
\end{array}
\right. $

\strut

(5.12) $\ \ \ =\left\{ 
\begin{array}{ll}
\left( 2N\right) ^{\frac{n}{2}}\cdot c_{\frac{n}{2}} & \text{if }n\text{ is
even} \\ 
&  \\ 
0 & \text{if }n\text{ is odd,}
\end{array}
\right. $

\strut

where $c_{k}=\frac{1}{k+1}\left( 
\begin{array}{l}
2k \\ 
\,\,k
\end{array}
\right) $ is the $k$-th Catalan number, for all $k\in \mathbb{N}.$ Remember
that

\strut

$\ \ \ \ \ \ \ \ \ \left| NC(k)\right| =\left| NC_{2}(2k)\right| =c_{k},$
for all $k\in \mathbb{N}.$

\strut

Therefore, by (5.9) and (5.12), we can compute the moments and cumulants of
the generating operator $T$ of $\left( W^{*}(G),tr\right) $ ;

\strut

$\ \ \ \ \ \ \ \ \ tr\left( T^{n}\right) =\left\{ 
\begin{array}{ll}
\left( 2N\right) ^{\frac{n}{2}}\cdot c_{\frac{n}{2}} & \text{if }n\text{ is
even} \\ 
&  \\ 
0 & \text{if }n\text{ is odd,}
\end{array}
\right. $

and

$\ \ \ \ \ \ \ \ \ k_{n}\left( T,...,T\right) =\left\{ 
\begin{array}{lll}
2N &  & \text{if }n=2 \\ 
&  &  \\ 
0 &  & \text{otherwise.}
\end{array}
\right. $

\strut
\end{example}

\strut

\begin{example}
Let $N\in \mathbb{N}$ and let $G$ be the circulant graph with

\strut

$\ \ \ \ \ \ \ \ \ \ \ \ \ V(G)=\{v_{1},...,v_{N}\}$

and

$\ \ \ \ \ \ \ \ \ \ \ \ \ E(G)=\{e_{1},...,e_{N}\}$

with

\strut

$\ \ e_{j}=v_{j}e_{j}v_{j+1}$, for $j=1,...,N-1,$ and $%
e_{N}=v_{N}e_{N}v_{1}. $

\strut

Define the generating operator $T=\sum_{j=1}^{N}\left(
L_{e_{j}}+L_{e_{j}}^{*}\right) $ of the graph $W^{*}$-algebra $W^{*}(G).$ In
this case, we can get the diagonal subalgebra $D_{G}$ of $W^{*}(G),$ as a
von Neumann algebra which is isomorphic to $\Delta _{N},$ where $\Delta _{N}$
is a subalgebra of the matricial algebra $M_{N}(\mathbb{C}).$ Define the
canonical conditional expectation $E:W^{*}(G)\rightarrow D_{G}.$ Then we can
compute the trivial $n$-th $D_{G}$-valued cumulant of the operator $T,$ by
regarding it as a $D_{G}$-valued random variable in the graph $W^{*}$%
-probability space $\left( W^{*}(G),E\right) $ over $D_{G}=\Delta _{N}.$
Notice that each block $L_{e_{j}}+L_{e_{j}}^{*}$'s are free from each other
over $D_{G}$ in $\left( W^{*}(G),E\right) ,$ by the diagram-distinctness of $%
e_{j}$'s, for $j=1,...,N.$

\strut

Fix $n\in \mathbb{N}.$ Then

\strut

$\ k_{n}\left( \underset{n\text{-times}}{\underbrace{T,.......,T}}\right)
=k_{n}\left( \sum_{j=1}^{N}\left( L_{e_{j}}+L_{e_{j}}^{*}\right)
,...,\sum_{j=1}^{N}\left( L_{e_{j}}+L_{e_{j}}^{*}\right) \right) $

\strut

$\ \ \ \ \ \ \ \ \ \ \ =\sum_{j=1}^{N}k_{n}\left(
(L_{e_{j}}+L_{e_{j}}^{*}),...,(L_{e_{j}}+L_{e_{j}}^{*})\right) $

\strut

by the mutual $D_{G}$-freeness of $\{L_{e_{j}},L_{e_{j}}^{*}\}$'s, for $%
j=1,...,N$

\strut

(5.13) $\ =\sum_{j=1}^{N}\underset{(u_{1},...,u_{n})\in \{1,*\}}{\sum }%
k_{n}\left( L_{e_{j}}^{u_{1}},...,L_{e_{j}}^{u_{n}}\right) .$

\strut

Recall that, by Section 2.2, we can get that

\strut

(5.14) $\ \ k_{n}\left( L_{e_{j}}^{u_{1}},...,L_{e_{j}}^{u_{n}}\right) =\mu
_{e_{j},...,e_{j}}^{u_{1},...,u_{n}}\cdot \Pr oj\left(
L_{e_{j}}^{u_{1}}...L_{e_{j}}^{u_{n}}\right) ,$

\strut

where $\mu _{e_{j},...,e_{j}}^{u_{1},...,u_{n}}=\underset{\pi \in
C_{e_{j},...,e_{j}}^{u_{1},...,u_{n}}}{\sum }\mu (\pi ,1_{n}).$

\strut

Observe that since $e_{j}$'s are non-loop edges, $e_{j}^{k}\notin \mathbb{F}%
^{+}(G),$ for all $k\in \mathbb{N}\,\setminus \,\{1\}$, for $j=1,...,N.$ In
other words, such $e_{j}^{k}$ is not admissible. So, if $(u_{1},...,u_{n})$
is not alternating, in the sense that $(u_{1},...,u_{n})=(1,*,...,1,*)$ or $%
(*,1,...,*,1),$ then $\Pr oj\left(
L_{e_{j}}^{u_{1}}...L_{e_{j}}^{u_{n}}\right) =0_{D_{G}}.$ For instance, $%
E\left( L_{e_{j}}^{*}L_{e_{j}}L_{e_{j}}^{*}\right) =0_{D_{G}}$ or $E\left(
L_{e_{j}}^{2}L_{e_{j}}^{*}L_{e_{j}}\right) =0_{D_{G}},$ by Section 2.1.
Therefore, the only nonvanishing case is either

\strut

$\ \ \ k_{n}\left(
L_{e_{j}},L_{e_{j}}^{*},...,L_{e_{j}},L_{e_{j}}^{*}\right) $ \ or \ $%
k_{n}\left( L_{e_{j}}^{*},L_{e_{j}},...,L_{e_{j}}^{*},L_{e_{j}}\right) ,$

\strut

where $n$ is even. Notice that

\strut

(5.15) $\ \ \ \ \ \ \ \mu _{e_{j},e_{j},...,e_{j},e_{j}}^{1,*,...,1,*}=\mu
_{e_{j},e_{j},...,e_{j},e_{j}}^{*,1,...,*,1},$

\strut

because $%
C_{e_{j},e_{j},...,e_{j},e_{j}}^{1,*,...,1,*}=C_{e_{j},e_{j},...,e_{j},e_{j}}^{*,1,...,*,1}, 
$ for all $j=1,...,N.$ Moreover, since $%
C_{e_{j},e_{j},...,e_{j},e_{j}}^{1,*,...,1,*}$ $=$ $%
C_{e_{k},e_{k},...,e_{k},e_{k}}^{1,*,...,1,*},$ for all $j\neq k$ in $%
\{1,...,N\},$

\strut

(5.16) $\ \ \ \ \ \ \ \mu _{e_{j},e_{j},...,e_{j},e_{j}}^{1,*,...,1,*}=\mu
_{e_{k},e_{k},...,e_{k},e_{k}}^{1,*,...,1,*},$

\strut

for all $j,k\in \{1,...,N\}.$ Let's denote $\mu
_{e_{j},e_{j},...,e_{j},e_{j}}^{1,*,...,1,*}$ by $\mu _{n},$ for all $%
j=1,...,N.$ Then, by (5.14), we have that

\strut

\strut (5.17)

$k_{n}\left( L_{e_{j}}^{u_{1}},...,L_{e_{j}}^{u_{n}}\right) =\left\{ 
\begin{array}{ll}
\mu _{n}L_{v_{j}} & \text{if }(u_{1},...,u_{n})=(1,*,...,1,*) \\ 
&  \\ 
\mu _{n}L_{v_{j+1}} & \text{if }(u_{1},...,u_{n})=(*,1,...,*,1) \\ 
&  \\ 
0_{D_{G}} & \text{otherwise,}
\end{array}
\right. $

\strut

for all $j=1,...,N,$ where $L_{v_{N+1}}$ means $L_{v_{1}}.$ So, by (5.13)
and (5.17), we can get that

\strut

$\ k_{n}\left( T,...,T\right) $

\strut

$\ =\sum_{j=1}^{N}\left( k_{n}\left(
L_{e_{j}},L_{e_{j}}^{*}...,L_{e_{j}},L_{e_{j}}^{*}\right) +k_{n}\left(
L_{e_{j}}^{*},L_{e_{j}},...,L_{e_{j}}^{*},L_{e_{j}}\right) \right) $

\strut \strut

$\ =\sum_{j=1}^{N}\left( \mu _{n}L_{v_{j}}+\mu _{n}L_{v_{j+1}}\right)
=\sum_{j=1}^{N}\mu _{n}\left( L_{v_{j}}+L_{v_{j+1}}\right) $

\strut

where $L_{v_{N+1}}$ means $L_{v_{1}}$, for all $n\in 2\mathbb{N}.$ Therefore,

\strut

$\ \ k_{n}\left( T,...,T\right) =\left\{ 
\begin{array}{ll}
\sum_{j=1}^{N}\mu _{n}\left( L_{v_{j}}+L_{v_{j+1}}\right) & \text{if }n\text{
is even} \\ 
&  \\ 
0_{D_{G}} & \text{if }n\text{ is odd.}
\end{array}
\right. $

\strut

$\ $(5.18)$\ \ \ \ \ \ \ \ =\left\{ 
\begin{array}{ll}
2\mu _{n}\cdot 1_{D_{G}} & \text{if }n\text{ is even} \\ 
&  \\ 
0_{D_{G}} & \text{if }n\text{ is odd.}
\end{array}
\right. $

\strut

Unfortunately, it is very hard to compute $\mu _{n},$ when $n\rightarrow
\infty .$ But we have to remark that if we have arbitrary graph $H$ and its
graph $W^{*}$-probability space $\left( W^{*}(H),F\right) $ over its
diagonal subalgebra $D_{H}$ and if $w\in loop^{c}(G),$ then

\strut

$\ \ \ \ \ \ \ \ \ \mu _{w,w,...,w,w}^{1,*,...,1,*}=\mu _{n}=\mu
_{w,w,...,w,w}^{*,1,...,*,1},$ for all $n\in 2\mathbb{N}.$

\strut

Now, let's compute the trivial $n$-th $D_{G}$-valued moment of $T.$ Notice
that since all odd trivial $D_{G}$-valued cumulants of $T$ vanish, all odd
trivial $D_{G}$-valued moments of $T$ vanish (See [11] and [14]). Thus it
suffices to compute the even trivial $D_{G}$-valued moments of $T.$ Assume
that $n\in 2\mathbb{N}.$ Then

\strut

(5.19)\ \ \ \ \ \ \ \ \ $\ E\left( T^{n}\right) =\underset{\pi \in NC_{E}(n)%
}{\sum }k_{\pi }\left( T,...,T\right) ,$

\strut

where $k_{\pi }\left( T,...,T\right) $ is the partition-dependent cumulant
of $T$ (See [16]) and

\strut

$\ \ \ \ NC_{E}(n)\overset{def}{=}\{\pi \in NC(n):V\in \pi \Leftrightarrow
\left| V\right| \in 2\mathbb{N}\}.$

\strut

By (5.18), we can get that $k_{n}(T,...,T)$ commutes with all elements in $%
W^{*}(G),$ because $1_{D_{G}}$ and $0_{D_{G}}$ commutes with $W^{*}(G)$ and $%
2\mu _{n}\in \mathbb{C},$ for all $n\in \mathbb{N}.$ So, the formula (5.19)
can be reformed by

\strut

$\ \ \ E(T^{n})=\underset{\pi \in NC_{E}(n)}{\sum }\,\left( \underset{V\in
\pi }{\prod }k_{\left| V\right| }(T,...,T)\right) $

\strut

$\ \ \ \ \ \ \ \ \ \ \ \ =\underset{\pi \in NC_{E}(n)}{\sum }\,\left( 
\underset{V\in \pi }{\prod }2\mu _{\left| V\right| }\cdot 1_{D_{G}}\right) $

\strut

(5.20)$\ \ \ \ \ =\left( \underset{\pi \in NC_{E}(n)}{\sum }\,\left( 
\underset{V\in \pi }{\prod }2\mu _{\left| V\right| }\right) \right) \cdot
1_{D_{G}},$

\strut

for all $n\in 2\mathbb{N}.$ Therefore, by (5.18) and (5.20), we have that if 
$T$ is the generating operator of the graph $W^{*}$-algebra of the circulant
graph $G$ with $N$-vertices, then

\strut \strut \strut

$\ \ \ E(T^{n})=\left\{ 
\begin{array}{ll}
\left( \underset{\pi \in NC_{E}(n)}{\sum }\,\left( \underset{V\in \pi }{%
\prod }2\mu _{\left| V\right| }\right) \right) \cdot 1_{D_{G}} & \text{if }n%
\text{ is even} \\ 
&  \\ 
0_{D_{G}} & \text{if }n\text{ is odd.}
\end{array}
\right. $

and

$\ \ \ \ \ \ \ \ \ k_{n}\left( \underset{n\text{-times}}{\underbrace{%
T,.....,T}}\right) =\left\{ 
\begin{array}{ll}
2\mu _{n}\cdot 1_{D_{G}} & \text{if }n\text{ is even} \\ 
&  \\ 
0_{D_{G}} & \text{if }n\text{ is odd.}
\end{array}
\right. $

\strut \strut $\ $
\end{example}

\strut \strut \strut

\strut \strut

\begin{quote}
\textbf{Reference}

\strut

\strut

{\small [1] \ \ A. Nica, R-transform in Free Probability, IHP course note,
available at www.math.uwaterloo.ca/\symbol{126}anica.}

{\small [2]\strut \ \ \ A. Nica and R. Speicher, R-diagonal Pair-A Common
Approach to Haar Unitaries and Circular Elements, (1995), www
.mast.queensu.ca/\symbol{126}speicher.\strut }

{\small [3] \ }$\ ${\small B. Solel, You can see the arrows in a Quiver
Operator Algebras, (2000), preprint}

{\small \strut [4] \ \ A. Nica, D. Shlyakhtenko and R. Speicher, R-cyclic
Families of Matrices in Free Probability, J. of Funct Anal, 188 (2002),
227-271.}

{\small [5] \ \ D. Shlyakhtenko, Some Applications of Freeness with
Amalgamation, J. Reine Angew. Math, 500 (1998), 191-212.\strut }

{\small [6] \ \ D.Voiculescu, K. Dykemma and A. Nica, Free Random Variables,
CRM Monograph Series Vol 1 (1992).\strut }

{\small [7] \ \ D. Voiculescu, Operations on Certain Non-commuting
Operator-Valued Random Variables, Ast\'{e}risque, 232 (1995), 243-275.\strut 
}

{\small [10]\ D. Shlyakhtenko, A-Valued Semicircular Systems, J. of Funct
Anal, 166 (1999), 1-47.\strut }

{\small [10]\ D.W. Kribs and M.T. Jury, Ideal Structure in Free Semigroupoid
Algebras from Directed Graphs, preprint}

{\small [10]\ D.W. Kribs and S.C. Power, Free Semigroupoid Algebras, preprint%
}

{\small [11]\ I. Cho, Amalgamated Boxed Convolution and Amalgamated
R-transform Theory, (2002), preprint.}

{\small [12] I. Cho, The Tower of Amalgamated Noncommutative Probability
Spaces, (2002), Preprint.}

{\small [13] I. Cho, Free Perturbed R-transform Theory, (2003), Preprint.}

{\small [14]\ I. Cho, Compatibility of a Noncommutative Probability Space
and a Noncommutative Probability Space with Amalgamation, (2003), Preprint. }

{\small [15] I. Cho, Graph }$W^{*}${\small -Probability Spaces Over the
Diagonal Subalgebras, (2004), Preprint.}

{\small [16] R. Speicher, Combinatorial Theory of the Free Product with
Amalgamation and Operator-Valued Free Probability Theory, AMS Mem, Vol 132 ,
Num 627 , (1998).}

{\small [17] R. Speicher, Combinatorics of Free Probability Theory IHP
course note, available at www.mast.queensu.ca/\symbol{126}speicher.\strut }

{\small [18] T. Bates and D. Pask, Flow Equivalence of Graph Algebras,
(2004), Preprint.}
\end{quote}

\end{document}